\documentclass[review,hidelinks,onefignum,onetabnum]{siamart220329}

\usepackage{lipsum}
\usepackage{amsfonts}
\usepackage{graphicx}
\usepackage{epstopdf}
\usepackage{algorithmic}
\ifpdf
  \DeclareGraphicsExtensions{.eps,.pdf,.png,.jpg}
\else
  \DeclareGraphicsExtensions{.eps}
\fi

\newsiamremark{remark}{Remark}
\newsiamremark{hypothesis}{Hypothesis}
\crefname{hypothesis}{Hypothesis}{Hypotheses}
\newsiamthm{claim}{Claim}

\headers{SAV-based entropy-dissipative schemes for kinetic equations}{S. Zhang, J. Shen, and J. Hu}

\title{SAV-based entropy-dissipative schemes for a class of kinetic equations\thanks{Submitted to the editors DATE.
\funding{The work of SZ and JH is supported by DOE grant DE-SC0023164. JH is also partially supported by NSF grant DMS-2409858 and AFOSR grant FA9550-21-1-0358. The work of JS is supported by National Natural Science Foundation of China under grant No. 12371409.}}}

\author{Shiheng Zhang\thanks{Department of Applied Mathematics, University of Washington, Seattle, WA 98195, USA
  (\email{shzhang3@uw.edu}).}
\and Jie Shen\thanks{School of Mathematical Science,
Eastern Institute of Technology, Ningbo, Zhejiang 315200, P. R. China (\email{jshen@eitech.edu.cn}).}
\and Jingwei Hu\thanks{Corresponding author. Department of Applied Mathematics, University of Washington, Seattle, WA 98195, USA (\email{hujw@uw.edu}).}}

\usepackage{amsopn}

\ifpdf
\hypersetup{
  pdftitle={SAV-based entropy-dissipative schemes for a class of kinetic equations},
  pdfauthor={Shiheng Zhang, Jie Shen, and Jingwei Hu}
}
\fi

\usepackage{subcaption}
\newcommand{\rd}{\,\mathrm{d}}

\begin{document}\nolinenumbers

\maketitle

\begin{abstract}
We introduce novel entropy-dissipative numerical schemes for a class of kinetic equations, leveraging the recently introduced scalar auxiliary variable (SAV) approach. Both first and second order schemes are constructed. Since the positivity of the solution is closely related to entropy, we also propose positivity-preserving versions of these schemes to ensure robustness, which include a scheme specially designed for the Boltzmann equation and a more general scheme using Lagrange multipliers. The accuracy and provable entropy-dissipation properties of the proposed schemes are validated for both the Boltzmann equation and the Landau equation through extensive numerical examples.

\end{abstract}



\section{Introduction}
\label{sec:intro}

We are interested in structure-preserving discretizations to the prototype kinetic equation\footnote{Many kinetic equations also depend on the spatial variable $x$. In this work, we restrict ourselves to the spatially homogeneous case, that is, $f(t,x,v)$ is assumed to be homogeneous in the $x$ direction.} given by
\begin{equation} \label{kinetic}
\partial_tf =Q(f), 
\end{equation}
where $f=f(t,v)$ is the probability density function of time $t\geq 0$ and velocity $v\in \mathbb{R}^d$, 
and $Q(f)$ is the collision operator modeling particle collisions. Depending on the application, $Q(f)$ can take various forms, but they generally satisfy the following common properties:
\begin{itemize}
\item conservation of mass, momentum, and energy:
\begin{equation} \label{conservation}
\int_{\mathbb{R}^d}Q(f)\rd{v}=\int_{\mathbb{R}^d}Q(f)v\rd{v}=\int_{\mathbb{R}^d}Q(f)v^2\rd{v}=0,
\end{equation}
\item Boltzmann's {\it H-theorem}:
\begin{equation} \label{entropy}
\int_{\mathbb{R}^d} Q(f)\log f\rd{v}\leq 0,
\end{equation}
and the equality sign holds if and only if $f$ becomes the {\it Maxwellian}:
\begin{equation}
M[f]=\frac{\rho}{(2\pi T)^{d/2}}\exp \left( -\frac{|v-u|^2}{2T}\right),
\end{equation}
where the density $\rho$, bulk velocity $u$, and temperature $T$ are defined as
\begin{equation}
\rho = \int_{\mathbb{R}^d} f\rd{v}, \quad u=\frac{1}{\rho} \int_{\mathbb{R}^d} fv\rd{v}, \quad T=\frac{1}{d\rho} \int_{\mathbb{R}^d} f|v-u|^2\rd{v}.
\end{equation}
\end{itemize}
Using these properties, one can easily show that the solution to equation \eqref{kinetic} satisfies
\begin{equation}
\frac{\rd }{\rd t} \int_{\mathbb{R}^d} f (1,v,|v|^2)^T\rd{v}=\int_{\mathbb{R}^d} \partial_t f(1,v,|v|^2)^T \rd{v}=\int_{\mathbb{R}^d} Q(f)(1,v,|v|^2)^T \rd{v}= 0,
\end{equation}
hence $\rho$, $u$, and $T$ remain as constant; furthermore,
\begin{equation} \label{eq:entropy}
\frac{\rd }{\rd t} \int_{\mathbb{R}^d} f\log f\rd{v}=\int_{\mathbb{R}^d} \partial_t f(\log f+1)\rd{v}=\int_{\mathbb{R}^d} Q(f)(\log f+1)\rd{v}\leq 0,
\end{equation}
hence the entropy $\int_{\mathbb{R}^d} f\log f\rd{v}$ decays over time.

A large class of kinetic equations falls into the form of \eqref{kinetic}, including the Boltzmann equation \cite{Cercignani}, the Landau equation \cite{Landau37}, and their simplified versions such as the BGK model \cite{BGK54}, the ES-BGK model \cite{Holway66}, and the kinetic Fokker-Planck model \cite{Villani02}, to name a few. In particular, the Boltzmann collision operator that describes neutral particle collisions is given by 
\begin{equation} \label{Boltzcol}
Q_B(f)(v)=\int_{\mathbb{R}^d}\int_{S^{d-1}}B(v-v_*,\sigma)[f(v')f(v_*')-f(v)f(v_*)]\rd{\sigma}\rd{v_*},
\end{equation}
where $\sigma$ is a vector varying over the unit sphere $S^{d-1}$ ($d\geq 2$), $(v,v_*)$ and $(v',v_*')$ are the velocity pairs before and after a collision, related by
\begin{equation}
v'=\frac{v+v_*}{2}+\frac{|v-v_*|}{2}\sigma, \quad v_*'=\frac{v+v_*}{2}-\frac{|v-v_*|}{2}\sigma,
\end{equation}
and the collision kernel $B$ assumes the form
\begin{equation} \label{kernel}
B(v-v_*,\sigma)=C_B|v-v_*|^{\gamma}b(\cos \theta), \quad \cos \theta=\frac{\sigma \cdot (v-v_*)}{|v-v_*|}, \quad -d< \gamma \leq 1.
\end{equation}
On the other hand, the Landau collision operator that describes charged particle collisions is given by 
\begin{equation} \label{Landcol}
Q_L(f)(v) = \nabla_v \cdot \int_{\mathbb{R}^d} A(v-v_*)[f(v_*)\nabla_v f(v) - f(v)\nabla_{v_*} f(v_*)] \rd{v_*},
\end{equation}
where the collision kernel $A$ is a $d\times d$ positive semi-definite matrix given by
\begin{equation}
A(v-v_*)=C_L|v-v_*|^{\gamma}(|v-v_*|^2I_d-(v-v_*)\otimes (v-v_*)), \quad -d\leq \gamma \leq 1.
\end{equation}

When numerically solving the equation \eqref{kinetic}, it is desirable to have a time discretization scheme that captures the entropy-decay structure \eqref{eq:entropy} at the discrete level. We first note that this can be easily achieved via the backward Euler scheme:
\begin{equation}
\frac{f^{n+1}-f^n}{\Delta t}=Q(f^{n+1}),
\end{equation}
which implies
\begin{equation}
\int_{\mathbb{R}^d} f^{n+1}\log f^{n+1}\rd{v}-\int_{\mathbb{R}^d} f^n\log f^{n+1}\rd{v}=\Delta t \int_{\mathbb{R}^d}Q(f^{n+1})\log f^{n+1}\rd{v}\leq 0,
\end{equation}
hence
\begin{equation}
\int_{\mathbb{R}^d} f^{n+1}\log f^{n+1}\rd{v}\leq \int_{\mathbb{R}^d} f^n\log f^{n+1}\rd{v} \leq \int_{\mathbb{R}^d} f^n\log f^{n}\rd{v},
\end{equation}
where the last inequality is due to Jensen's inequality. However, the backward Euler scheme can be extremely difficult to implement for the Boltzmann equation and Landau equation due to their nonlocal and nonlinear collision operators. On the other hand, the forward Euler scheme applied to \eqref{kinetic}, although simple to implement, generally cannot preserve the entropy-decay structure\footnote{There are some exceptions: e.g., the forward Euler for the Boltzmann equation with  Maxwell collision kernel (i.e., $\gamma=0$ in \eqref{kernel}) is readily entropy-dissipative, as proved in \cite{Villani98}.}.

Motivated by the above discussion, our goal in this work is to design entropy-dissipative schemes for general kinetic equations \eqref{kinetic} that are as easily implementable as the forward Euler. This presents a significant challenge, primarily due to the inherent complexity of the collision operator which precludes the application of some well-known strategies to achieve energy stability (a property similar to entropy dissipation), such as convex splitting \cite{elliott1993global,eyre1998unconditionally,shen2012second}. Recently, the scalar auxiliary variable (SAV) method \cite{shen2018scalar} has emerged as a successful approach in various domains, notably within the contexts of gradient flows \cite{shen2018convergence, shen2019new, cheng2018multiple, liu2023novel}. The method has subsequently been extended to solve a variety of complex problems, such as the Schrödinger equation \cite{antoine2021scalar}, Navier-Stokes equations \cite{lin2019numerical}, and Wasserstein gradient flows \cite{zhang2024structure}. Given the simplicity and efficiency of the SAV method, it has also prompted extensive research beyond solving PDEs, including optimization \cite{liu2023efficient, zhang2023element} and machine learning \cite{ma2024efficient, zhang2024energy}.

This paper proposes a novel strategy to develop entropy-dissipative schemes for kinetic equations \eqref{kinetic} leveraging the SAV approach. Both first and second order schemes are constructed. Since the positivity of the solution is closely related to entropy (as $\log f$ is not even defined for negative $f$), we also propose positivity-preserving versions of these schemes to ensure robustness, which include a scheme specially designed for the Boltzmann equation and a more general scheme using Lagrange multipliers \cite{ito2008lagrange,harker1990finite,facchinei2003finite,bergounioux1999primal}. This latter approach has recently been applied to parabolic problems to address challenges in preserving positivity and mass conservation, as demonstrated in \cite{van2019positivity, cheng2022new}.



Finally, we mention that our main focus of this work is regarding the time discretization. To simplify the presentation, we keep the velocity variable $v$ as continuous in the following discussion. However, it should be understood that a discretization method for the collision operator is employed and we assume that this method satisfies the following properties:
\begin{itemize}
\item it assumes a truncation in the velocity domain: $v\in \Omega$ instead of $\mathbb{R}^d$;
\item it preserves mass, momentum and energy in $\Omega$ as in \eqref{conservation}; 
\item it satisfies the {\it H-theorem} in $\Omega$ as in \eqref{entropy}. 
\end{itemize}
There exist such methods in the literature, for example, the discrete velocity methods for the Boltzmann collision operator \cite{DP14} and the Landau collision operator \cite{DL94} satisfy the above assumptions. The Fourier spectral methods for the Boltzmann collision operator \cite{PR00,MP06,GHHH17} and the Landau collision operator \cite{PRT00} coupled with the entropy fix \cite{CHKL22} can also satisfy the above assumptions, except for momentum and energy conservation.

The rest of this paper is organized as follows. Section \ref{sec:entropySAV} introduces the entropy-dissipative SAV schemes. These schemes can achieve entropy dissipation but the positivity of the solution is not guaranteed. Subsequently, Section \ref{sec:boltz} presents a positivity-preserving and entropy-dissipative SAV scheme, tailored specifically for the Boltzmann equation. Section \ref{sec:PPentropySAV} introduces positivity-preserving and entropy-dissipative SAV schemes applicable to a broader class of kinetic equations.  Section \ref{sec:num} provides a summary of the numerical schemes proposed in this paper along with their corresponding properties, and  presents extensive numerical examples for the Boltzmann equation and Landau equation to demonstrate the performance of the proposed schemes. Some concluding remarks are given in Section \ref{sec:con}.


\section{Entropy-dissipative SAV schemes}
\label{sec:entropySAV}



We define $H=\int_{\Omega} f\log f\rd v + C$, where $C$ is a constant chosen such that $H\geq H_{\min} >0$ for all $t\geq 0$. This is a reasonable assumption since the entropy $\int_{\Omega} f\log f\rd v$ decreases over time and attains its minimum when the system reaches the steady-state. We then introduce a scalar auxiliary variable (SAV) $r$: $$r = \sqrt{H},$$ and rewrite the equation \eqref{kinetic} equivalently as
\begin{align}
    \partial_tf & = \frac{r}{\sqrt{H}}Q(f), \label{sav-cont1}\\
    \frac{\rd{r}}{\rd{t}}&=\frac{1}{2\sqrt{H}}\int_{\Omega}\log f\partial_t f\rd{v}, \label{sav-cont2}
\end{align}
where the second equation is derived by differentiating \( r \) with respect to $t$ and applying conservation of mass.

\subsection{A first order scheme}


A first order numerical scheme, \texttt{SAV-1st}, for \eqref{sav-cont1}-\eqref{sav-cont2} can be obtained as follows:
\begin{align}
    \frac{f^{n+1}-f^n}{\Delta t}  &= \frac{r^{n+1}}{\sqrt{H^n}} Q(f^n) \label{sav-landau-1},\\
    \frac{r^{n+1}-r^n}{\Delta t} &= \frac{1}{2\sqrt{H^n}} \int_{\Omega}\log f^n \frac{f^{n+1}-f^{n}}{\Delta t} \rd v,\label{sav-landau-2}
\end{align}
where $H^n=H(f^n)=\int_{\Omega} f^n\log f^n \rd v + C>0$, and $r^0= \sqrt{H^0}$. Note that in general $r^n \neq \sqrt{H^n}$ for $n>0$. This scheme is easy to implement. One can determine \( r^{n+1} \) by plugging \eqref{sav-landau-1} into \eqref{sav-landau-2}, and then compute \( f^{n+1} \) using \eqref{sav-landau-1}.

\begin{theorem}\label{Thm-SAV-1st}The scheme \eqref{sav-landau-1}-\eqref{sav-landau-2}, \texttt{SAV-1st}, satisfies the following properties: for all time steps $n\geq 0$ and step size $\Delta t>0$,
\begin{itemize}
\item conservation of mass, momentum and energy: $$\int_{\Omega} f^{n+1}(1,v,|v|^2)^T\rd{v}=\int_{\Omega} f^{n}(1,v,|v|^2)^T \rd{v};$$
\item modified entropy dissipation: $$\tilde{H}^{n+1} - \tilde{H}^{n} =-(r^{n+1}-r^n)^2+ \frac{\Delta t(r^{n+1})^2}{H^n} \int_{\Omega}Q(f^n) 
    \log f^n \rd v \leq 0,$$
where $\tilde{H}^{n} = (r^n)^2$. 
\end{itemize}
\end{theorem}
\begin{proof}

  Conservation of mass, momentum, and energy can be easily shown by multiplying \eqref{sav-landau-1} by $(1,v,|v|^2)^T$ and integrating in $v$, and using that $$\int_{\Omega}Q(f^n)(1,v,|v|^2)^T\rd{v}=0.$$
  
    To show the entropy dissipation, we can multiply \eqref{sav-landau-1} by $\frac{r^{n+1}}{\sqrt{H^n}}\log f^n$ and integrate in $v$, and multiply \eqref{sav-landau-2} by $2r^{n+1}$ to obtain
    \begin{align*}
    \frac{r^{n+1}}{\sqrt{H^n}}\int_{\Omega}\log f^n \frac{f^{n+1}-f^{n}}{\Delta t} \rd v  &= \frac{(r^{n+1})^2}{H^n} \int_{\Omega}Q(f^n) 
    \log f^n \rd v,\\
    \frac{2r^{n+1}(r^{n+1}-r^n)}{\Delta t} &= \frac{r^{n+1}}{\sqrt{H^n}} \int_{\Omega}\log f^n \frac{f^{n+1}-f^{n}}{\Delta t} \rd v.
\end{align*}
Combining these two equations together, we have
\begin{equation*}
    \frac{2r^{n+1}(r^{n+1}-r^n)}{\Delta t} = \frac{(r^{n+1})^2}{H^n} \int_{\Omega}Q(f^n) 
    \log f^n \rd v.
\end{equation*}
Noting the identity
\begin{equation}\label{identity1}
    2r^{n+1}(r^{n+1}-r^n)  = (r^{n+1})^2 - (r^n)^2 + (r^{n+1}-r^n)^2,
\end{equation}
we then have 
\begin{equation*}
    (r^{n+1})^2-(r^{n})^2 + (r^{n+1}-r^n)^2 = \frac{\Delta t(r^{n+1})^2}{H^n} \int_{\Omega}Q(f^n) 
    \log f^n \rd v \leq 0.
\end{equation*}
\end{proof}
\subsection{A second order scheme}
A second order numerical scheme, \texttt{SAV-2nd}, for \eqref{sav-cont1}-\eqref{sav-cont2} can be constructed by extending the second order Backward Differentiation Formula (BDF) as follows:
\begin{align}
    \frac{3f^{n+1}-4f^n+f^{n-1}}{2\Delta t}  &= \frac{r^{n+1}}{\sqrt{H^{n+1,*}}} Q(f^{n+1,*}),\label{sav2-landau-1}\\
    \frac{3r^{n+1}-4r^n+r^{n-1}}{2\Delta t} &= \frac{1}{2\sqrt{H^{n+1,*}}} \int_{\Omega}\log f^{n+1,*}\frac{3f^{n+1}-4f^{n}+f^{n-1}}{2\Delta t} \rd v,\label{sav2-landau-2}
\end{align}
where $r^0 = \sqrt{H^0}$ and $r^1$, $f^1$ are obtained by the first order scheme \texttt{SAV-1st}. Further, $H^{n+1,*}$ represents $H(f^{n+1,*})$, where $f^{n+1,*}$ is an explicit approximation of $f(t^{n+1})$ with an order of accuracy $\mathcal{O}(\Delta t^2)$. For example, the Adams-Bashforth extrapolation \cite{butcher2016numerical} can be employed for this purpose:
\begin{equation*}
     f^{n+1,*} = 2 f^n - f^{n-1} .
\end{equation*}
\begin{theorem}\label{Thm-SAV-2nd}The scheme \eqref{sav2-landau-1}-\eqref{sav2-landau-2}, \texttt{SAV-2nd}, satisfies the following properties: for all time steps $n\geq 0$ and step size $\Delta t>0$,
\begin{itemize}
\item conservation of mass, momentum and energy: $$\int_{\Omega} f^{n+1}(1,v,|v|^2)^T\rd{v}=\int_{\Omega} f^{n}(1,v,|v|^2)^T \rd{v};$$ 
\item modified entropy dissipation: 
\begin{equation*}
    \begin{aligned}
        \tilde{H}^{n+1} - \tilde{H}^{n} &=-\frac{1}{2}(r^{n+1}-2r^n+r^{n-1})^2\\&+ \frac{\Delta t(r^{n+1})^2}{H^{n+1,*}} \int_{\Omega}Q(f^{n+1,*}) 
    \log f^{n+1,*} \rd v \leq 0,
    \end{aligned}
\end{equation*}
    where $\tilde{H}^{n} = \frac{1}{2}(r^n)^2 + \frac{1}{2}(2r^n-r^{n-1})^2$. 
\end{itemize}
\end{theorem}
\begin{proof}
  Conservation of mass, momentum, and energy can be easily shown by multiplying  \eqref{sav2-landau-1} by $(1,v,|v|^2)^T$ and integrating in $v$, and using that $$\int_{\Omega}Q(f^{n+1,*})(1,v,|v|^2)^T\rd{v}=0.$$
  
 To show the entropy dissipation, we multiply \eqref{sav2-landau-1} by $\frac{r^{n+1}}{\sqrt{H^{n+1,*}}}\log f^{n+1,*}$ and integrate in $v$, multiply \eqref{sav2-landau-2} by $2r^{n+1}$ to obtain
 \begin{align*}
    \frac{2r^n(3r^{n+1}-4r^n+r^{n-1})}{2\Delta t} = \frac{(r^{n+1})^2}{H^{n+1,*}} \int_{\Omega} Q(f^{n+1,*})\log f^{n+1,*}\rd{v}\leq 0.
 \end{align*}
Using the identity
\begin{equation}
\begin{aligned}
     2 r^{n+1}(3r^{n+1}-4r^n+r^{n-1})&=(r^{n+1})^2+(2r^{n+1}-r^n)^2+(r^{n+1}-2r^n+r^{n-1})^2\\
     &-(r^n)^2-(2r^n-r^{n-1})^2,
\end{aligned}
\end{equation}
we obtain the desired inequality.
\end{proof}

\begin{remark}
A second order scheme based on the Crank-Nicolson method can also be constructed:
\begin{align}
     &\frac{{f}^{n+1}-f^n}{\Delta t} =  \frac{r^{n+1} + r^n}{2\sqrt{H^{n+1/2,*}}} Q(f^{n+1/2,*}),\\
     &\frac{{r}^{n+1}-r^n}{\Delta t} = \frac{1}{2\sqrt{H^{n+1/2,*}}} \int_{\Omega}\log f^{n+1/2,*}\frac{{f}^{n+1}-f^n}{\Delta t} \rd v,
\end{align}
where $f^{n+1/2,*}$ represents any explicit $\mathcal{O}(\Delta t^2)$ approximation to $f(t^{n+1/2})$. This scheme satisfies the same properties as 
Theorem \ref{Thm-SAV-2nd}, except that the entropy dissipation has a different form:
\begin{equation*}
    \begin{aligned}
        (r^{n+1})^2 - (r^{n})^2 &= \frac{\Delta t(r^{n+1}+r^{n})^2}{4H^{n+1/2,*}} \int_{\Omega}Q(f^{n+1/2,*}) 
    \log f^{n+1/2,*} \rd v \leq 0.
    \end{aligned}
\end{equation*}
\end{remark}
\section{A positivity-preserving and entropy-dissipative SAV scheme for the Boltzmann equation}
\label{sec:boltz}
Since $f$ is the probability density function, so it should satisfy $f\ge 0$. 
In this section, we use the Boltzmann collision operator \eqref{Boltzcol} as an example to show how to modify the scheme \texttt{SAV-1st} in the last section to preserve positivity. First note that we can write $Q_B(f)$ as
\begin{equation}
Q_B(f) = Q_B^+(f) - Q_B^-(f)f, 
\end{equation}
with
\begin{align} 
Q_B^+(f)(v)&=\int_{\mathbb{R}^d}\int_{S^{d-1}}B(v-v_*,\sigma)f(v')f(v_*')\rd{\sigma}\rd{v_*}\geq 0,\\
Q_B^-(f)(v)&=\int_{\mathbb{R}^d}\int_{S^{d-1}}B(v-v_*,\sigma)f(v_*)\rd{\sigma}\rd{v_*}\geq 0.
\end{align}
We propose a first order stabilized scheme with positivity-preserving property, \\ \texttt{SAV-1st-P-B}, as follows
\begin{align}
    \frac{f^{n+1}-f^n}{\Delta t}  &=  \frac{{r}^{n+1}}{\sqrt{H^n}} Q_B(f^n) +  \beta f^n -  \beta f^{n+1},\label{psav-boltz-1}\\
    \frac{{r}^{n+1}-r^n}{\Delta t} &= \frac{1}{2\sqrt{H^n}} \int_{\Omega}\log f^n \frac{f^{n+1}-f^{n}}{\Delta t} \rd v,\label{psav-boltz-2}
\end{align}
where the constant $\beta$ is a stabilizing constant, chosen such that \\$\beta\geq \frac{r^0}{\sqrt{H_{\min}}}\max Q_B^-(f)\geq 0$.

\begin{theorem}The scheme \eqref{psav-boltz-1}-\eqref{psav-boltz-2}, \texttt{SAV-1st-P-B}, satisfies the following \\properties: for all time step $n\geq 0$ and step size $\Delta t>0$,
\begin{itemize}
\item the solution $f^n\geq 0$, provided initially $f^0 \geq 0$;
\item conservation of mass, momentum and energy: $$\int_{\Omega} f^{n+1}(1,v,|v|^2)^T\rd{v}=\int_{\Omega} f^{n}(1,v,|v|^2)^T \rd{v};$$
\item modified entropy dissipation: $$\tilde{H}^{n+1} - \tilde{H}^{n} =-(r^{n+1}-r^n)^2+ \frac{\Delta t(r^{n+1})^2}{(1+\beta \Delta t)H^n} \int_{\Omega}Q_B(f^n) 
    \log f^n \rd v \leq 0,$$
where $\tilde{H}^{n} = (r^n)^2$. 
\end{itemize}
\end{theorem}

\begin{proof}

     To show the positivity of $f^n$, we first rewrite the equation 
       \eqref{psav-boltz-1} as
       \begin{equation}
       \begin{aligned}
            \frac{f^{n+1}-f^n}{\Delta t}  &=  \frac{{r}^{n+1}}{\sqrt{H^n}} \left(Q_B^+(f^n) - Q_B^-(f^n) f^n\right) +  \beta f^n -  \beta f^{n+1}\\
            &=\frac{{r}^{n+1}}{\sqrt{H^n}} Q_B^+(f^n) + \left(\beta-\frac{{r}^{n+1}}{\sqrt{H^n}}Q_B^-(f^n)\right) f^n -  \beta f^{n+1},
       \end{aligned}
       \end{equation}
       hence
       \begin{equation}\label{psav-boltz-p}
            (1+{\Delta t}\beta){f^{n+1}} = 
             \Delta t\frac{{r}^{n+1}}{\sqrt{H^n}} Q_B^+(f^n) + \left(1 +  \Delta t\left(\beta-\frac{{r}^{n+1}}{\sqrt{H^n}}Q_B^-(f^n)\right)\right) f^n.
       \end{equation}
Rewriting \eqref{psav-boltz-1} alternatively as
  \begin{equation}\label{psav-boltz-3}
       \left(\frac{1}{\Delta t} + \beta\right)\left(f^{n+1}-f^n\right) = \frac{{r}^{n+1}}{\sqrt{H^n}} Q_B(f^n),
  \end{equation}   
  and plugging \eqref{psav-boltz-3} into \eqref{psav-boltz-2}, we can obtain:
  \begin{equation}
      r^{n+1} = \left(1-\frac{\Delta t}{2{H^n}(1+\Delta t\beta)} \int_{\Omega}  Q_B(f^n)\log f^n \rd v  \right)^{-1} r^n.
  \end{equation}
This implies that $0\leq r^{n+1} \leq r^n$.
By the choice of $\beta$, we guarantee that 
$$\beta \geq \frac{r^0}{\sqrt{H_{\min}}}\max Q_B^-(f)\geq \frac{{r}^{n+1}}{\sqrt{H^n}}Q_B^-(f^n).$$
Using this in \eqref{psav-boltz-p}, we see that  $f^{n+1}\geq 0$ if $f^n\geq 0$.
     
 Conservation of mass, momentum, and energy is immediate by multiplying equation \eqref{psav-boltz-1} by $(1,v,|v|^2)^T$ and integrating in $v$.

To show the entropy dissipation, we can multiply \eqref{psav-boltz-3} by $\frac{r^{n+1}}{\sqrt{H^n}}\log f^n$ and integrate in $v$, and multiply  \eqref{psav-boltz-2} by $2r^{n+1}$ to obtain
\begin{align*}
    \left(\frac{1}{\Delta t} + \beta\right)\frac{r^{n+1}}{\sqrt{H^n}}\int_{\Omega}\log f^n \left(f^{n+1}-f^{n}\right) \rd v  &= \frac{(r^{n+1})^2}{H^n} \int_{\Omega}Q_B(f^n) 
    \log f^n \rd v,\\
    \frac{2r^{n+1}(r^{n+1}-r^n)}{\Delta t} &= \frac{r^{n+1}}{\sqrt{H^n}} \int_{\Omega}\log f^n \frac{f^{n+1}-f^{n}}{\Delta t} \rd v.
\end{align*}
Combining the two equations together, we have
\begin{equation*}
    2r^{n+1}(r^{n+1}-r^n) = \left(\frac{1}{\Delta t} + \beta \right)^{-1}\frac{(r^{n+1})^2}{H^n} \int_{\Omega}Q_B(f^n) 
    \log f^n \rd v.
\end{equation*}
Using the identity \eqref{identity1} together with $\beta > 0$, we have 
\begin{equation*}
    (r^{n+1})^2-(r^{n})^2 + (r^{n+1}-r^n)^2 = \left(\frac{1}{\Delta t} + \beta \right)^{-1}\frac{(r^{n+1})^2}{H^n} \int_{\Omega}Q_B(f^n) 
    \log f^n \rd v \leq 0.
\end{equation*}

\end{proof}


\section{Positivity-preserving schemes for general kinetic equations}
\label{sec:PPentropySAV}

The \\scheme introduced in Section~\ref{sec:boltz} is designed for the Boltzmann equation and is limited to first order. In this section, we construct positivity-preserving and entropy-dissipative SAV schemes that work for general kinetic equations \eqref{kinetic}, leveraging the optimization techniques. We will achieve this in two steps. The first version focuses on restoring the positivity without mass conservation; and the second version can achieve both positivity and mass conservation.

\subsection{Positivity-preserving schemes without mass conservation}
\subsubsection{A first order scheme}

To guarantee that \(f\) remains positive, we introduce a Lagrange multiplier function, \(\lambda(t,v)\), and consider the extended system with the Karush–Kuhn–Tucker (KKT) conditions:
\begin{align}
   & \partial_tf - Q(f) = \lambda,\\
   & \lambda \geq 0, \quad  f\geq 0, \quad  \lambda f = 0.
\end{align}
A first order operator splitting scheme \cite{cheng2022new} with SAV and Lagrange multiplier, \\ \texttt{SAV-1st-L}, is given as follows:

\textbf{Step 1} (prediction): solve $\tilde{f}^{n+1}$ from
\begin{align}
     &\frac{\tilde{f}^{n+1}-f^n}{\Delta t} =  \frac{{r}^{n+1}}{\sqrt{H^n}} Q(f^n) \label{SAV-1st-L1},\\
     &\frac{{r}^{n+1}-r^n}{\Delta t} = \frac{1}{2\sqrt{H^n}} \int_{\Omega}\log f^n \frac{\tilde{f}^{n+1}-f^{n}}{\Delta t} \rd v\label{SAV-1st-L2}.
\end{align}

\textbf{Step 2} (correction): solve $\left(f^{n+1}, \lambda^{n+1}\right)$ from
\begin{align}
     &\frac{{f}^{n+1}(v)-\tilde{f}^{n+1}(v)}{\Delta t} = \lambda^{n+1}(v)\label{SAV-1st-L3},\\
    &\lambda^{n+1}(v) \geq 0, \quad f^{n+1}(v)\geq 0, \quad \lambda^{n+1}(v) f^{n+1}(v) = 0.\label{SAV-1st-L4}
\end{align}

The equations in the prediction step is the same as the scheme \texttt{SAV-1st}, hence can be solved directly and enjoy the entropy-decay property. A notable characteristic of the equations in the correction step is their solvability on a point-wise basis, as described below:
\begin{equation}\label{cut-off1}
    \left(f^{n+1}(v), \lambda^{n+1}(v)\right)=\left\{\begin{array}{cl}
\left(\tilde{f}^{n+1}(v), 0\right) & \text { if } \tilde{f}^{n+1}(v)\geq 0 \\
\left(0,-\frac{\tilde{f}^{n+1}(v)}{\Delta t}\right) & \text { otherwise }
\end{array}, \quad \forall \ v \in \Omega.\right.
\end{equation}

\begin{theorem} \label{Thm-SAV-1st-L}
The scheme \eqref{SAV-1st-L1}-\eqref{SAV-1st-L4}, \texttt{SAV-1st-L}, satisfies the following properties: for all time steps $n\geq 0$ and step size $\Delta t>0$,
\begin{itemize}
\item the solution $f^n\geq 0$; 
\item modified entropy dissipation: $$\tilde{H}^{n+1} - \tilde{H}^{n} =-(r^{n+1}-r^n)^2+ \frac{\Delta t(r^{n+1})^2}{H^n} \int_{\Omega}Q(f^n) 
    \log f^n \rd v \leq 0,$$
where $\tilde{H}^{n} = (r^n)^2$. 
\end{itemize}
\end{theorem}
\begin{proof}

It is clear that $f^n\geq 0$ for all $n$. As for the modified entropy dissipation, it can be established  following a similar approach to that detailed in Theorem \ref{Thm-SAV-1st}.

\end{proof}
\begin{remark}
    We can also require that the solution is bounded away from 0 for a prescribed $\epsilon$ by substituting the optimality condition \eqref{SAV-1st-L4} with
    \begin{equation}\label{postive kkt}
        \lambda^{n+1}(v) \geq 0,\quad  f^{n+1}(v)\geq \epsilon, \quad \lambda^{n+1}(v)(f^{n+1}(v)- \epsilon)=0.
   \end{equation}
\end{remark}

\subsubsection{A second order scheme}
A second order scheme with SAV and Lagrange multiplier, \texttt{SAV-2nd-L}, can also be constructed with the second order BDF and Adams-Bashforth extrapolation:

\textbf{Step 1} (prediction): solve $\tilde{f}^{n+1}$ from
\begin{align}
      &\frac{3\tilde{f}^{n+1}-4f^n + f^{n-1}}{2\Delta t} =  \frac{r^{n+1}}{\sqrt{H^{n+1,*}}} Q(f^{n+1,*}) ,\label{SAV-2nd-L1}\\
     &\frac{3{r}^{n+1}-4r^n+r^{n-1}}{2\Delta t} = \frac{1}{2\sqrt{H^{n+1,*}}} \int_{\Omega}\log f^{n+1,*}\frac{3\tilde{f}^{n+1}-4f^{n}+f^{n-1}}{2\Delta t} \rd v\label{SAV-2nd-L2}.
\end{align}

\textbf{Step 2} (correction): solve $\left(f^{n+1}, \lambda^{n+1}\right)$ from
\begin{align}
      &\frac{3{f}^{n+1}(v)-3\tilde{f}^{n+1}(v)}{2\Delta t} = \lambda^{n+1}(v)\label{SAV-2nd-L3},\\
    &\lambda^{n+1}(v) \geq 0, \quad f^{n+1}(v)\geq 0, \quad \lambda^{n+1}(v) f^{n+1}(v) = 0.\label{SAV-2nd-L4}
\end{align}
Note that the Adams-Bashforth extrapolation can not preserve positivity, we need to modify it with
\begin{equation}
   f^{n+1,*} = 
\begin{cases} 
2 f^n - f^{n-1}, & \text{if } f^n \geq f^{n-1}, \\
\frac{1}{2 / f^n - 1 / f^{n-1}}, & \text{otherwise}.
\end{cases}
\end{equation}
The steps \eqref{SAV-2nd-L3} to \eqref{SAV-2nd-L4} can be solved in a similar way as in the first order case:
\begin{equation}\label{cut-off2}
\left(f^{n+1}(v), \lambda^{n+1}(v)\right)=\left\{\begin{array}{cl}
\left(\tilde{f}^{n+1}(v), 0\right) & \text { if } \tilde{f}^{n+1}(v)\geq 0 \\
\left(0,-\frac{3\tilde{f}^{n+1}(v)}{2\Delta t}\right) & \text { otherwise }
\end{array}, \quad \forall \ v \in \Omega.\right.
\end{equation}
\begin{theorem}\label{Thm-SAV-2nd-L}The scheme \eqref{SAV-2nd-L1}-\eqref{SAV-2nd-L4}, \texttt{SAV-2nd-L}, satisfies the following properties: for all time steps $n\geq 0$ and step size $\Delta t>0$,
\begin{itemize}
\item the solution $f^n\geq 0$;
\item modified entropy dissipation: 
      \begin{equation*}
       \begin{aligned}
            \tilde{H}^{n+1} - \tilde{H}^{n} &=-\frac{1}{2}(r^{n+1}-2r^n+r^{n-1})^2\\&+ \frac{\Delta t(r^{n+1})^2}{H^{n+1,*} }\int_{\Omega}Q(f^{n+1,*}) 
    \log f^{n+1,*} \rd v \leq 0,
       \end{aligned}
       \end{equation*}
    where $\tilde{H}^{n} = \frac{1}{2}(r^n)^2 + \frac{1}{2}(2r^n-r^{n-1})^2$. 
\end{itemize}
\end{theorem}
\begin{proof}
It is clear that $f^n\geq 0$ for all $n$.
  The proof for the modified entropy dissipation follows from a similar approach to that detailed in Theorem \ref{Thm-SAV-2nd}.
\end{proof}

\subsection{Positivity-preserving schemes with mass conservation}

The correction steps \eqref{cut-off1} and \eqref{cut-off2} are equivalent to the cut-off strategy \cite{lu2013cutoff}, which similarly encounter issues with mass conservation. For instance, we observe from \eqref{SAV-1st-L3} that
\begin{equation*}
\int_{\Omega} f^{n+1}(v) \rd v - \int_{\Omega} \tilde{f}^{n+1}(v) \rd v = \int_{\Omega} \Delta t \lambda^{n+1}(v) \rd v \geq 0,
\end{equation*}
which suggests an increase in mass. 

To enforce the mass conservation, we introduce an additional Lagrange multiplier, \(\xi^{n+1}\), which is independent of the velocity variable. This new multiplier aims to ensure mass conservation during the correction step.
A first order scheme that preserves positivity and mass conservation, \texttt{SAV-1st-LM}, is given as follows: 

\textbf{Step 1} (prediction): solve $\tilde{f}^{n+1}$ from
\begin{align}
     &\frac{\tilde{f}^{n+1}-f^n}{\Delta t} =  \frac{r^{n+1}}{\sqrt{H^n}} Q(f^n) \label{SAV-1st-LM1},\\
     &\frac{r^{n+1}-r^n}{\Delta t} = \frac{1}{2\sqrt{H^n}} \int_{\Omega}\log f^n \frac{\tilde{f}^{n+1}-f^{n}}{\Delta t} \rd v\label{SAV-1st-LM2}.
\end{align}

\textbf{Step 2} (correction): solve $\left(f^{n+1}, \lambda^{n+1},\xi^{n+1}\right)$ from
\begin{align}
      &\frac{{f}^{n+1}(v)-\tilde{f}^{n+1}(v)}{\Delta t} = \lambda^{n+1}(v) + \xi^{n+1}\label{SAV-1st-LM3},\\
    &\lambda^{n+1}(v) \geq 0,\quad  f^{n+1}(v)\geq 0, \quad \lambda^{n+1}(v) f^{n+1}(v) = 0\label{SAV-1st-LM4},\\
    & \int_{\Omega} {f}^{n+1}(v) \rd v =\int_{\Omega} {f}^{n}(v)  \rd v\label{SAV-1st-LM5}.
\end{align}

In order to solve the correction step, we rewrite \eqref{SAV-1st-LM3} in the following equivalent form
\begin{equation}\label{SAV-1st-LM3-1}
    \frac{f^{n+1}(v)-\left(\tilde{f}^{n+1}(v)+\Delta t \xi^{n+1}\right)}{\Delta t}=\lambda^{n+1}(v) .
\end{equation}
Hence, assuming $\xi^{n+1}$ is known, \eqref{SAV-1st-LM4} and \eqref{SAV-1st-LM3-1} can be solved point-wise similarly as in the previous subsection:
\begin{equation}\label{SAV-1st-LM4-1}
    \begin{aligned}
        \left(f^{n+1}(v), \lambda^{n+1}(v)\right) = \begin{cases}
            \left(\tilde{f}^{n+1}(v)+\Delta t \xi^{n+1}, 0\right) & \text{if } \tilde{f}^{n+1}(v)+\Delta t \xi^{n+1}\geq 0 \\
            \left(0,-\frac{\tilde{f}^{n+1}(v)+\Delta t \xi^{n+1}}{\Delta t}\right) & \text{otherwise}
        \end{cases}, \ \forall \ v \in \Omega.
    \end{aligned}
\end{equation}
It remains to determine $\xi^{n+1}$. We find from  \eqref{SAV-1st-LM5} and  \eqref{SAV-1st-LM3-1} that
$$
\int_{\Omega} \tilde{f}^{n+1}+\Delta t \xi^{n+1} \rd v = 
\int_{\Omega} f^n \rd v -
\int_{\Omega} \Delta t \lambda^{n+1} \rd v,
$$
which, thanks to  \eqref{SAV-1st-LM4-1}, can be rewritten as
\begin{equation}
    \int_{v \in \Omega \text { s.t. } 0<\tilde{f}^{n+1}(v)+\Delta t \xi^{n+1}} \tilde{f}^{n+1}+\Delta t \xi^{n+1} \rd v = 
\int_{\Omega} f^n \rd v.
\end{equation}
Hence, $\xi^{n+1}$ is a solution to the nonlinear algebraic equation
\begin{equation}
    F(\xi)=\int_{v \in \Omega \text { s.t. } 0<\tilde{f}^{n+1}(v)+\Delta t \xi}\tilde{f}^{n+1}+\Delta t \xi \rd v - 
\int_{\Omega} f^n \rd v = 0.
\end{equation}
Since $F^{\prime}(\xi)$ may not exist and is difficult to compute if it exists, instead of the Newton iteration, we can use the following secant method:
\begin{equation}
    \xi_{k+1}=\xi_k-\frac{F\left(\xi_k\right)\left(\xi_k-\xi_{k-1}\right)}{{F}\left(\xi_k\right)-{F}\left(\xi_{k-1}\right)} .
\end{equation}
Since $\xi^{n+1}$ is an approximation to zero, and it can be shown that $\xi^{n+1} \leq 0$ if we add \eqref{SAV-1st-LM3} to \eqref{SAV-1st-LM1} and take the integration, we can choose $\xi_0=0$ and $\xi_1=-O(\Delta t)$. Once $\xi^{n+1}$ is known, we can update $\left(f^{n+1}(v), \lambda^{n+1}(v)\right)$ with \eqref{SAV-1st-LM4-1}.

A second order scheme, \texttt{SAV-2nd-LM}, can be constructed as follows:

\textbf{Step 1} (prediction): solve $\tilde{f}^{n+1}$ from
\begin{align}
     &\frac{3\tilde{f}^{n+1}-4f^n + f^{n-1}}{2\Delta t} =  \frac{r^{n+1}}{\sqrt{H^{n+1,*}}} Q(f^{n+1,*}), \label{SAV-2nd-LM-1}\\
     &\frac{3r^{n+1}-4r^n+r^{n-1}}{2\Delta t} = \frac{1}{2\sqrt{H^{n+1,*}}} \int_{\Omega}\log f^{n+1,*}\frac{3\tilde{f}^{n+1}-4f^{n}+f^{n-1}}{2\Delta t} \rd v.\label{SAV-2nd-LM-2}
\end{align}

\textbf{Step 2} (correction): solve $\left(f^{n+1}, \lambda^{n+1}, \xi^{n+1}\right)$ from
\begin{align}
     &\frac{3{f}^{n+1}(v)-3\tilde{f}^{n+1}(v)}{2\Delta t} = \lambda^{n+1}(v)+ \xi^{n+1},\label{SAV-2nd-LM-3}\\
    &\lambda^{n+1}(v) \geq 0, \quad f^{n+1}(v)\geq 0, \quad \lambda^{n+1}(v) f^{n+1}(v) = 0,\label{SAV-2nd-LM-4}\\
    & \int_{\Omega} {f}^{n+1}(v) \rd v =\int_{\Omega} {f}^{n}(v)  \rd v,\label{SAV-2nd-LM-5}
\end{align}
where
\begin{equation}
   f^{n+1,*} = 
\begin{cases} 
2 f^n - f^{n-1}, & \text{if } f^n \geq f^{n-1}, \\
\frac{1}{2 / f^n - 1 / f^{n-1}}, & \text{otherwise}.
\end{cases}
\end{equation}

The first order scheme \eqref{SAV-1st-LM1}-\eqref{SAV-1st-LM5}, \texttt{SAV-1st-LM}, and the second order scheme \eqref{SAV-2nd-LM-1}-\eqref{SAV-2nd-LM-5}, \texttt{SAV-2nd-LM}, satisfy the same properties as in Theorem \ref{Thm-SAV-1st-L} and Theorem \ref{Thm-SAV-2nd-L}, respectively. In addition, they both conserve mass:
$$\int_{\Omega} f^{n+1}\rd{v}=\int_{\Omega} f^{n} \rd{v}.$$

\begin{remark}
The schemes introduced above only conserve mass. While conservation of momentum and energy could, in principle, be achieved by introducing additional Lagrange multipliers, this approach would lead to a coupled nonlinear system for the Lagrange multipliers, and complicate the solution process. Therefore, we do not pursue momentum and energy conservation here.
\end{remark}

\section{Numerical examples}
\label{sec:num}
In this section, we present several numerical results to demonstrate the properties of the proposed schemes.

For readers' convenience, we first provide a summary on the properties of the proposed numerical schemes in this paper. In Table \ref{properties}, ``Conservation" refers to conservation of mass, momentum, and energy (unless otherwise specified); ``modified entropy" could take different forms for different schemes; whenever a second order scheme is indicated, it means that the scheme satisfies the same properties as the first order scheme in the same row. 

\begin{table}[tbh]
\centering
\begin{tabular}{|p{3cm}||c|c|c|c|}
\hline
\centering  {Scheme} & {Conservation} & \centering\begin{tabular}[c]{@{}l@{}}\quad Modified \\ entropy decay \end{tabular} & {Positivity} & {Second Order} \\ \hline
\centering  \begin{tabular}[c]{@{}l@{}}\texttt{SAV-1st} \\ \eqref{sav-landau-1}-\eqref{sav-landau-2}\end{tabular}& yes & yes & no & \begin{tabular}[c]{@{}l@{}}\texttt{SAV-2nd} \\ \eqref{sav2-landau-1}-\eqref{sav2-landau-2}\end{tabular}\\ \hline
\centering \begin{tabular}[c]{@{}l@{}}\texttt{\quad SAV-1st-P-B} \\ \quad \quad \eqref{psav-boltz-1}-\eqref{psav-boltz-2}\\ only for Boltzmann \end{tabular} & yes & yes & yes & no \\ \hline
\centering  \begin{tabular}[c]{@{}l@{}}\texttt{SAV-1st-L} \\ \eqref{SAV-1st-L1}-\eqref{SAV-1st-L4}\end{tabular}& no & yes & yes &  \begin{tabular}[c]{@{}l@{}}\texttt{SAV-2nd-L} \\ \eqref{SAV-2nd-L1}-\eqref{SAV-2nd-L4}\end{tabular} \\ \hline
\centering  \begin{tabular}[c]{@{}l@{}}\texttt{SAV-1st-LM}\\ \eqref{SAV-1st-LM1}-\eqref{SAV-1st-LM5}\end{tabular}& only mass  & yes & yes & \begin{tabular}[c]{@{}l@{}}\texttt{SAV-2nd-LM}\\ \eqref{SAV-2nd-LM-1}-\eqref{SAV-2nd-LM-5}\end{tabular} \\ \hline
\end{tabular}
\caption{Properties of  the proposed schemes.}
\label{properties}
\end{table}

 In our numerical experiments below, we focus our attention on the Boltzmann collision operator \eqref{Boltzcol} and the Landau collision operator \eqref{Landcol}. For velocity domain discretization, we employ the Fourier spectral methods \cite{MP06,PRT00} for these operators. Although the Fourier spectral methods do not strictly satisfy the three conditions listed in Introduction, their high accuracy ensures that the error from velocity discretization is negligible compared to that from time discretization, allowing us to conduct a meaningful validation of the proposed schemes. For all of the following tests, we assume the two-dimensional velocity domain $\Omega=[-L,L)^2$ and use $N=64$ Fourier modes in each dimension. The model specific parameters are chosen as
\begin{itemize}
\item Boltzmann: collision kernel $ B=\frac{1}{2\pi}, L=(3\sqrt{2}+1)S/2;$
\item Landau: collision kernel
$A=\frac{1}{16}(|v-v_*|^2I_2-(v-v_*)\otimes (v-v_*)), L=2S.$
\end{itemize}
The value of $S$ will be specified in each test.

For all SAV-based schemes in the following, we set the constant $C = 10$ to ensure that $H(f)=\int_{\Omega} f\log f\rd v + C > 0$. 
For \texttt{sav-1st-LM} and \texttt{sav-2nd-LM}, our implementation ensures that \( f^{n} \geq \epsilon = 10^{-16} \), as described in equation \eqref{postive kkt}, to maintain well-defined logarithmic term, \( \log f^n \), throughout the computational process. To clarify the terminology used in the following numerical experiments, we recall the definitions of the entropies. The modified entropy for the SAV-based schemes is \((r^n)^2\) for first order schemes and \(\frac{1}{2}(r^{n})^2 + \frac{1}{2}(r^{n} - r^{n-1})^2\) for second order schemes. The actual entropy at time $t_n$ is \(H^n = H(f^n)\), where \(f^n\) is the numerical solution. The entropy of the analytical solution at time \(t_n\) is \(H(f(t_n))\), with \(f(t_n)\) being the analytical solution. 
\subsection{Test case 1: BKW solution}

The BKW solution is one of the few analytical solutions for the Boltzmann/Landau equation. When $d=2$, it is given by
\begin{equation}
f(t,v)=\frac{1}{2\pi K}\exp\left(-\frac{|v|^2}{2K}\right)\left(\frac{2K-1}{K}+\frac{1-K}{2K^2}|v|^2\right), \quad K=1-\exp(-t/8)/2.
\end{equation}
We take $t_0=0.5$ as the initial time and set $S=3.3$. Note that the same solution works for both Boltzmann and Landau for the aforementioned collision kernels. 

\subsubsection{Convergence tests}

\begin{figure}[htbp!]
\centering
\begin{subfigure}{0.48\textwidth}
    \includegraphics[width=\linewidth]{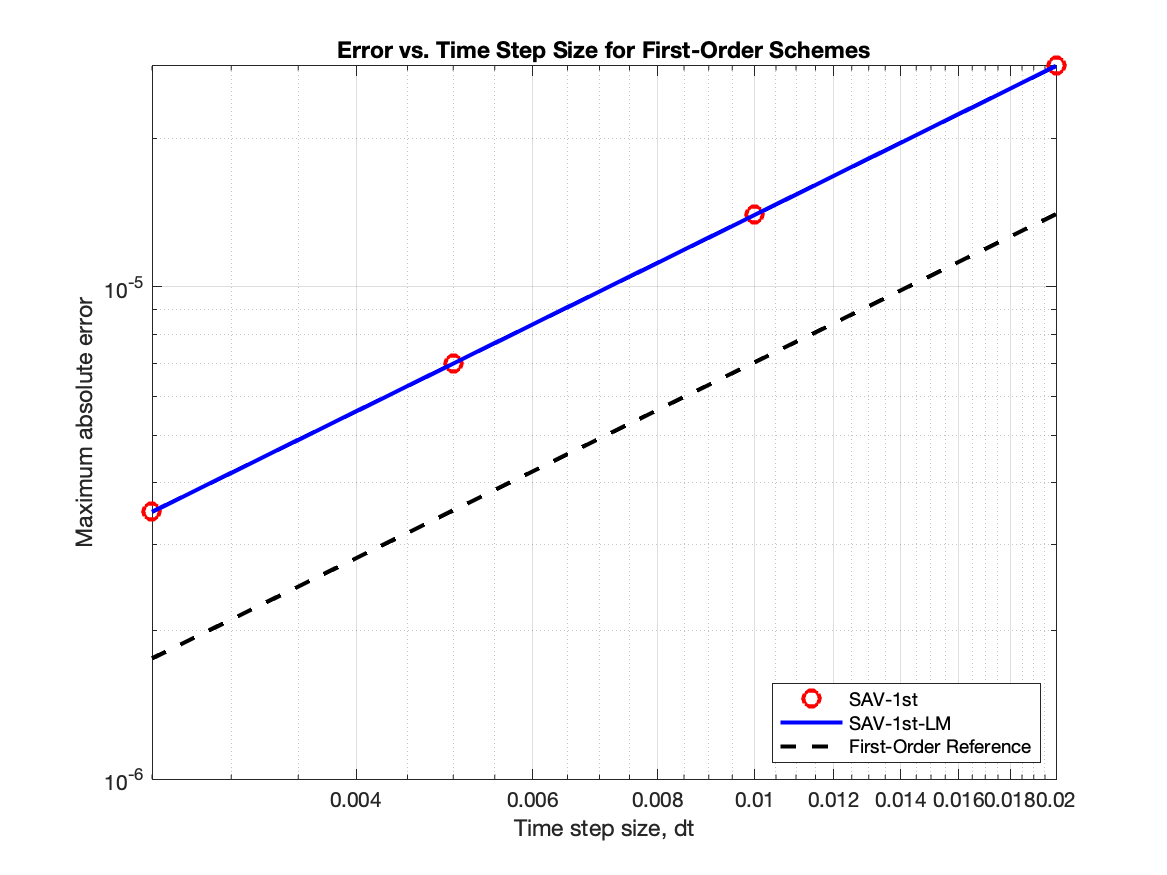}
    \caption{Boltzmann equation: error versus time step size for \texttt{SAV-1st} and \texttt{SAV-1st-LM}.}
    \label{fig:error_vs_dt_sav1st_comparison_boltz}
\end{subfigure}
\hfill
\begin{subfigure}{0.48\textwidth}
    \includegraphics[width=\linewidth]{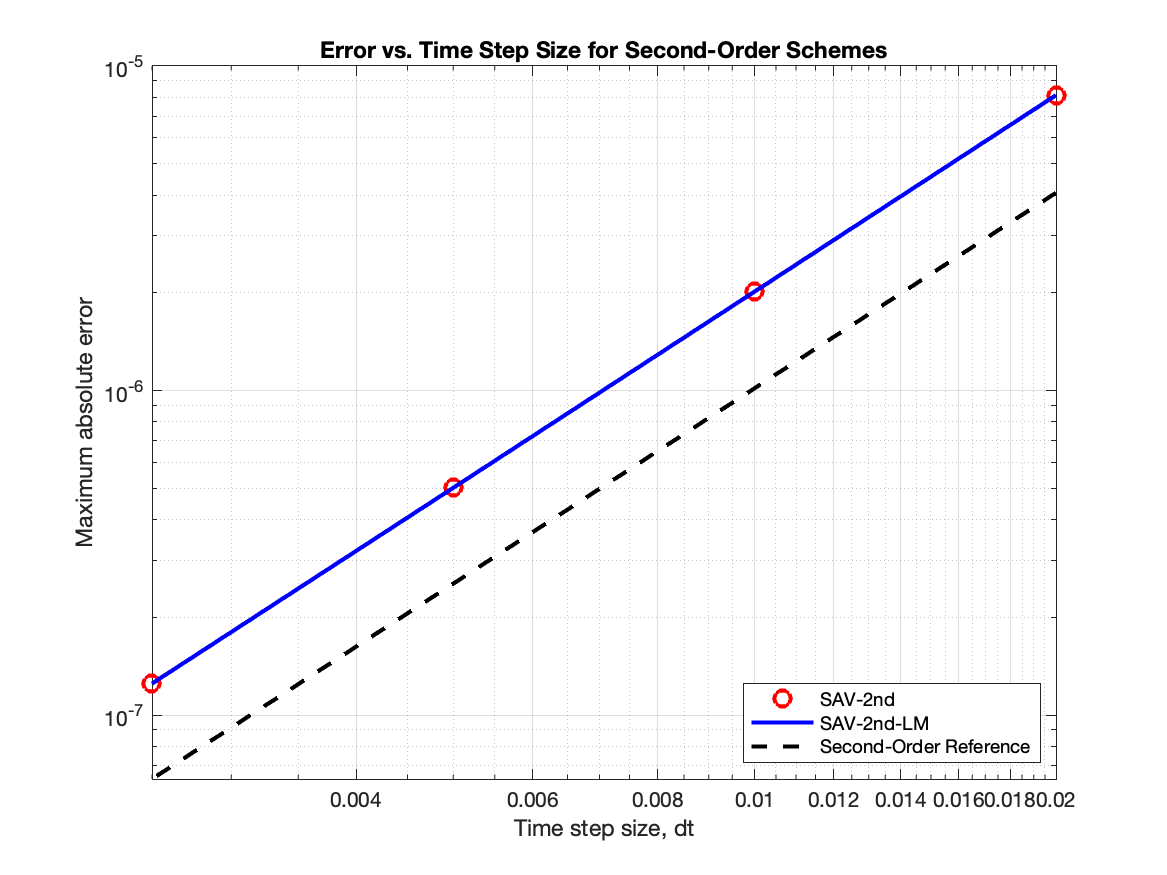}
    \caption{Boltzmann equation: error versus time step size for \texttt{SAV-2nd} and \texttt{SAV-2nd-LM}.}
    \label{fig:error_vs_dt_sav2nd_comparison_boltz}
\end{subfigure}
\hfill
\begin{subfigure}{0.48\textwidth}
    \includegraphics[width=\linewidth]{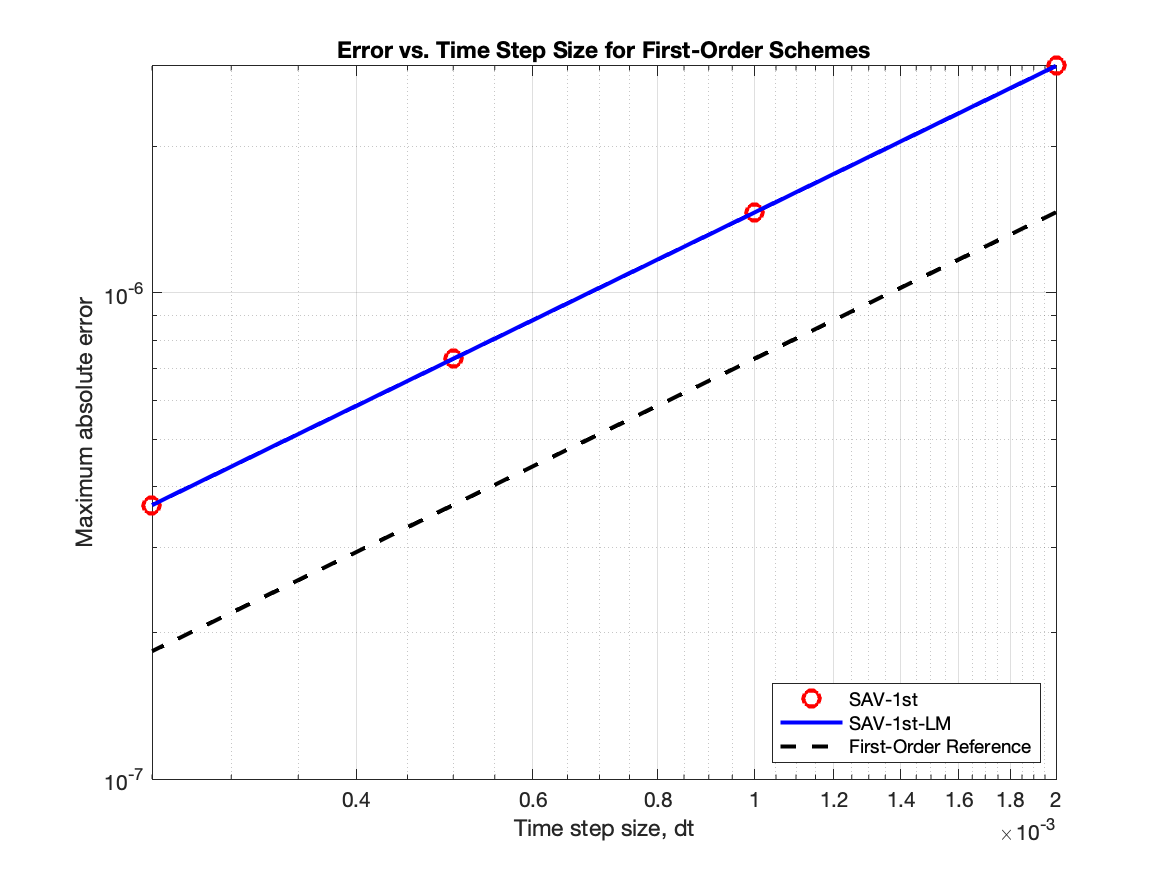}
    \caption{Landau equation: error versus smaller time step size for \texttt{SAV-1st} and \texttt{SAV-1st-LM}.}
    \label{fig:error_vs_dt_sav1st_comparison_landau}
\end{subfigure}
\hfill
\begin{subfigure}{0.48\textwidth}
    \includegraphics[width=\linewidth]{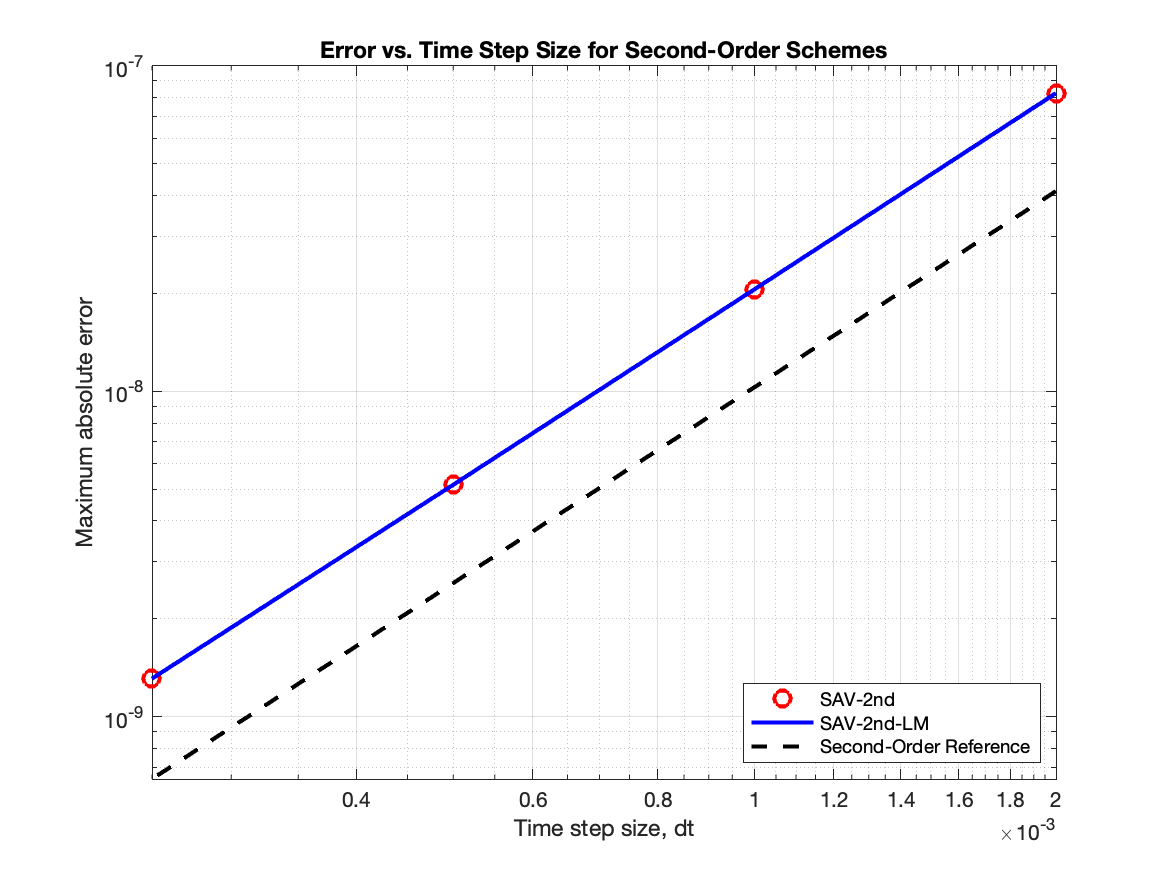}
    \caption{Landau equation: error versus smaller time step size for \texttt{SAV-2nd} and \texttt{SAV-2nd-LM}.}
    \label{fig:error_vs_dt_sav2nd_comparison_landau}
\end{subfigure}

\begin{subfigure}{0.48\textwidth}
    \includegraphics[width=\linewidth]{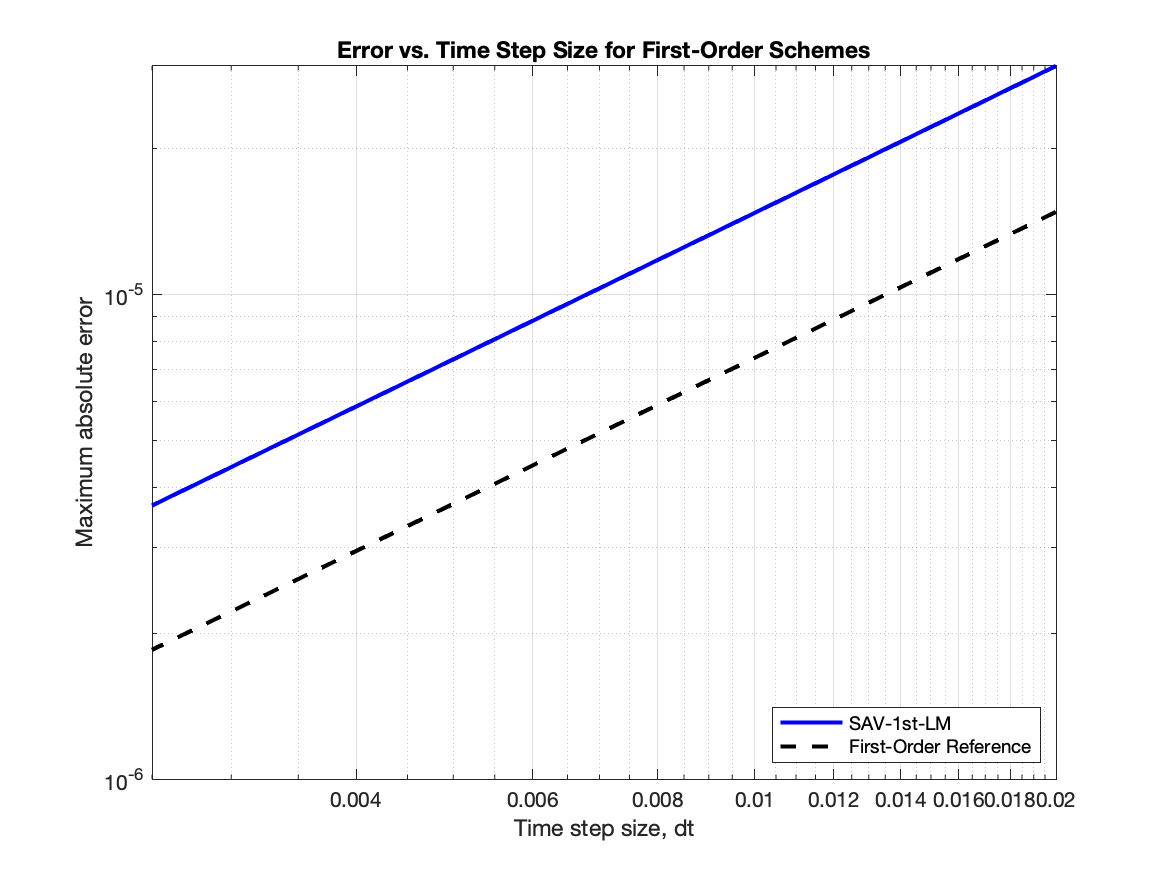}
    \caption{Landau equation: error versus larger time step size for \texttt{SAV-1st-LM}.}
\label{fig:error_vs_dt_sav1st_comparison_landau_lp}
\end{subfigure}
\hfill
\begin{subfigure}{0.48\textwidth}
    \includegraphics[width=\linewidth]{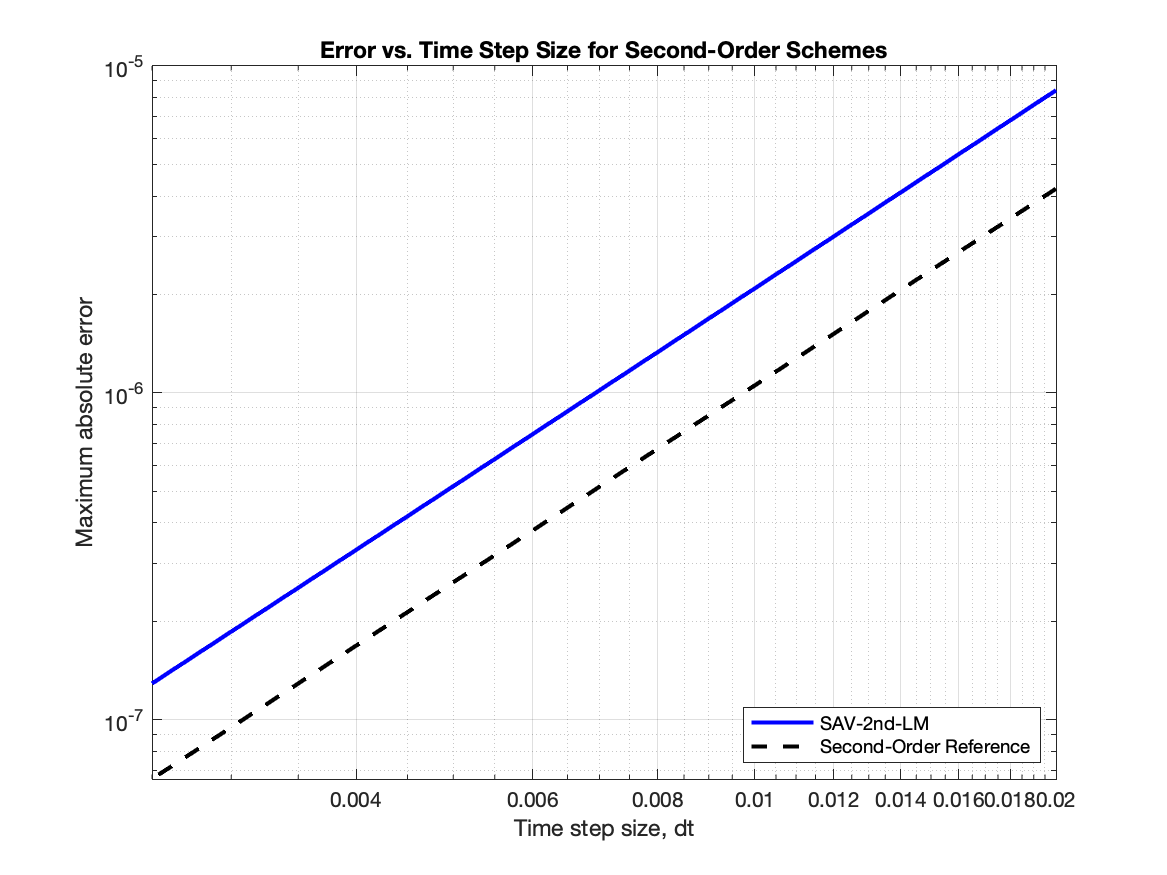}
    \caption{Landau equation: error versus larger time step size for \texttt{SAV-2nd-LM}.}
\label{fig:error_vs_dt_sav2nd_comparison_landau_lp}
\end{subfigure}
\caption{Convergence tests of \texttt{SAV-1st}, \texttt{SAV-2nd}, \texttt{SAV-1st-LM}, and \texttt{SAV-2nd-LM} for the Boltzmann equation and Landau equation.}

\label{fig:both_errors}
\end{figure}

We first perform a detailed error analysis to assess the convergence rates of \texttt{SAV-1st}, \texttt{SAV-2nd}, \texttt{SAV-1st-LM}, and \texttt{SAV-2nd-LM}. For the Boltzmann equation, time step sizes are set to $\Delta t = 0.02, 0.01, 0.005, 0.0025$, while for the Landau equation, much smaller steps of $\Delta t = 0.002, 0.001, 0.0005, 0.00025$ are used.
For each time step size, we run the solution up to the terminal time $t_{\text{end}}=0.6$ and record the maximum norm of the error between the numerical solution and the analytical solution. The results for the Boltzmann equation are shown in Figures \ref{fig:error_vs_dt_sav1st_comparison_boltz} and \ref{fig:error_vs_dt_sav2nd_comparison_boltz}, and those for the Landau equation are shown in Figures \ref{fig:error_vs_dt_sav1st_comparison_landau} and \ref{fig:error_vs_dt_sav2nd_comparison_landau}. The first and second order convergence are clearly observed.

It is worth noting that the time steps for the Landau equation in Figures \ref{fig:error_vs_dt_sav1st_comparison_landau} and \ref{fig:error_vs_dt_sav2nd_comparison_landau} are chosen to be small enough such that the positivity correction (hence mass conservation correction) is never triggered in \texttt{SAV-1st-LM} and \texttt{SAV-2nd-LM}. To further demonstrate the strength of these methods, we choose larger time step sizes of $\Delta t = 0.02, 0.01, 0.005, 0.0025$ and rerun the same test. The results are shown in Figures \ref{fig:error_vs_dt_sav1st_comparison_landau_lp} and \ref{fig:error_vs_dt_sav2nd_comparison_landau_lp}. In this case, \texttt{SAV-1st} and \texttt{SAV-2nd} result in negative values, causing the simulation to fail. In contrast, \texttt{SAV-1st-LM} and \texttt{SAV-2nd-LM} execute successfully and still demonstrate the expected convergence rates, highlighting their robustness under challenging conditions.

\subsubsection{Entropy evolution}

We then focus on the entropy evolution of the schemes \texttt{SAV-1st-LM} and \texttt{SAV-2nd-LM} for both the Boltzmann equation and Landau equation.

\begin{figure}[htbp!]
    \centering
    \begin{subfigure}[b]{0.45\textwidth}
        \includegraphics[width=\textwidth]{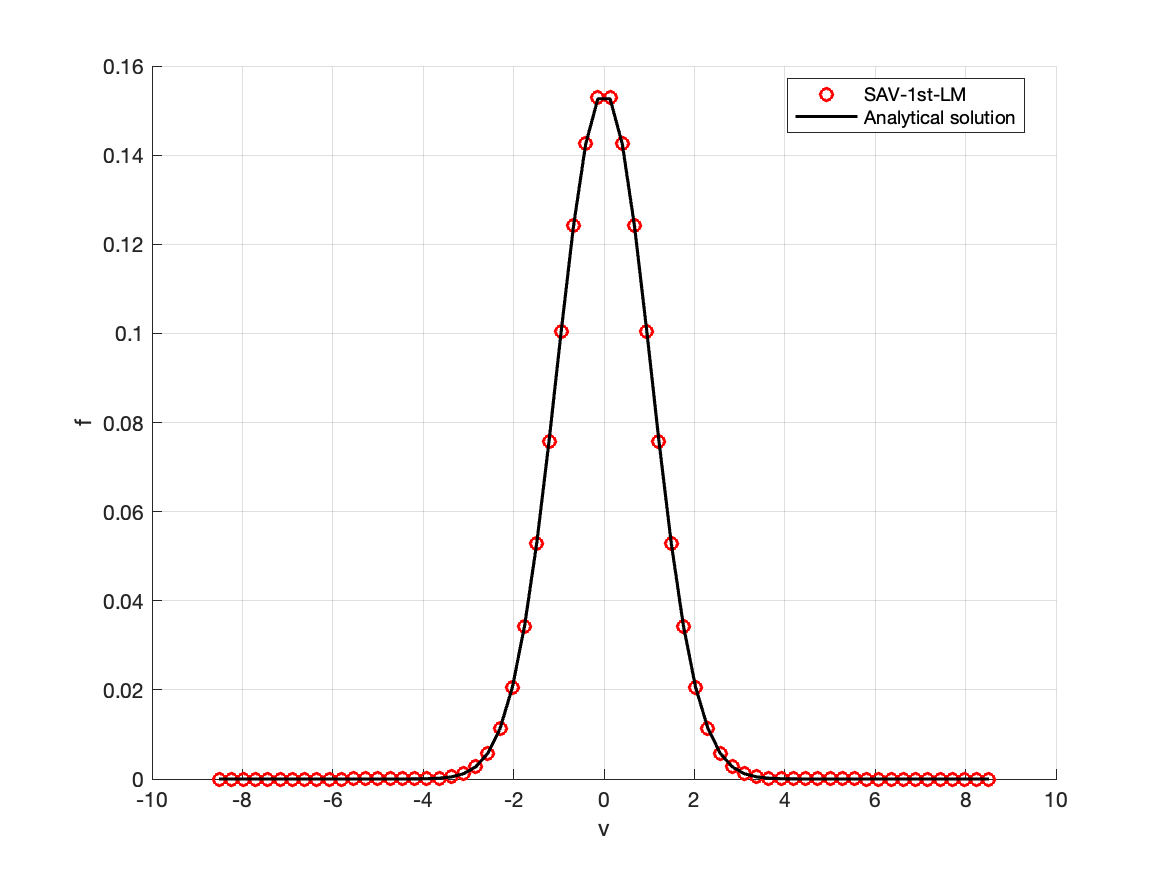}
        \caption{$f$ at $t_{\text{end}}=10.5$ with \( \Delta t = 0.2 \)}
    \end{subfigure}
    \hfill
    \begin{subfigure}[b]{0.45\textwidth}
        \includegraphics[width=\textwidth]{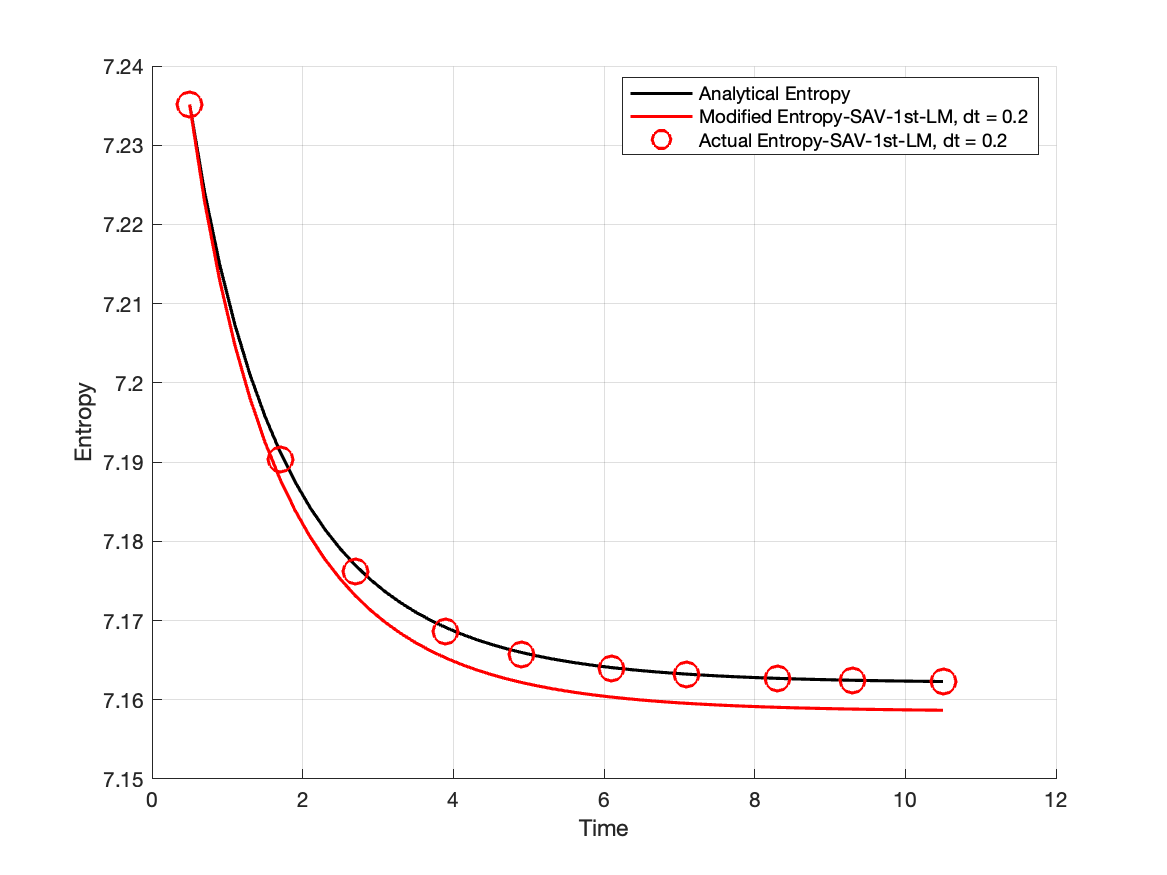}
        \caption{Entropy evolution with \(  \Delta t = 0.2 \)}
    \end{subfigure}
    \hfill
    \begin{subfigure}[b]{0.45\textwidth}
        \includegraphics[width=\textwidth]{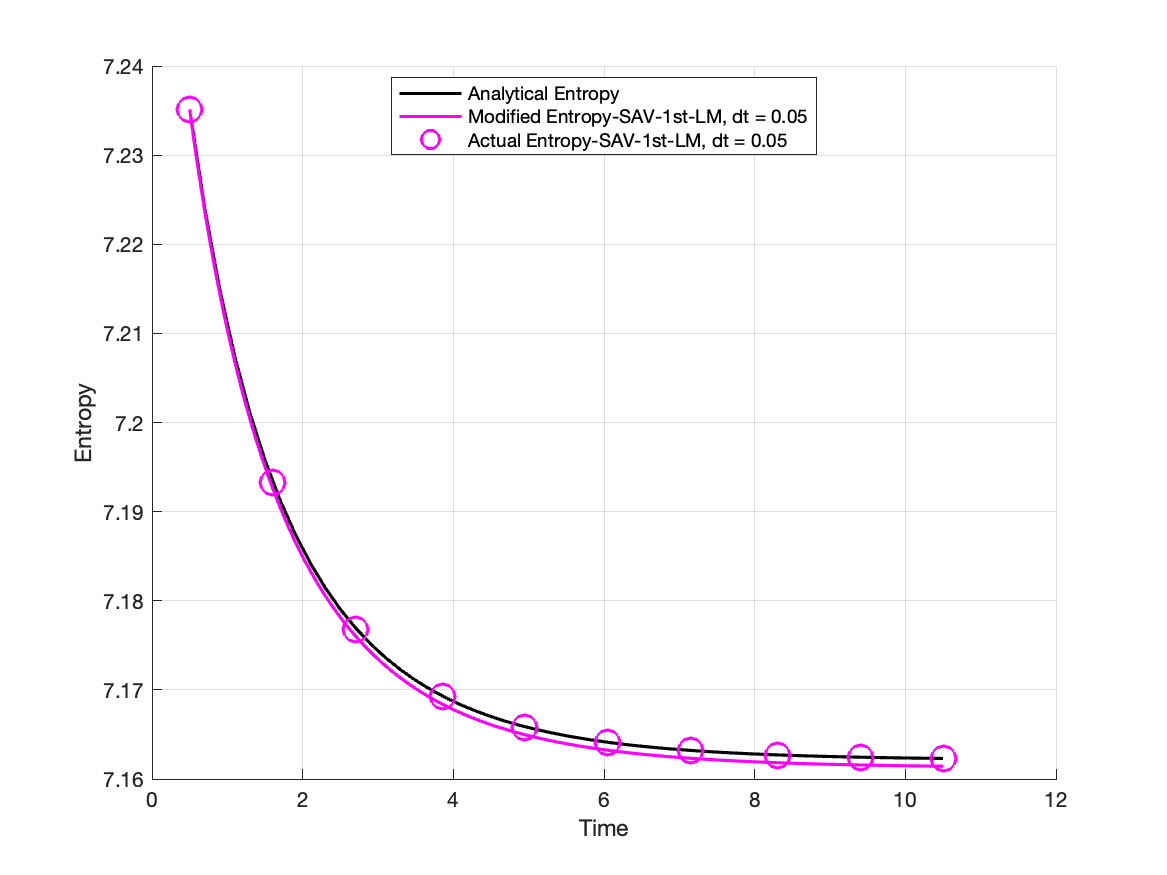}
        \caption{Entropy evolution with \( \Delta t = 0.05 \)}
    \end{subfigure}
    \hfill
    \begin{subfigure}[b]{0.45\textwidth}
        \includegraphics[width=\textwidth]{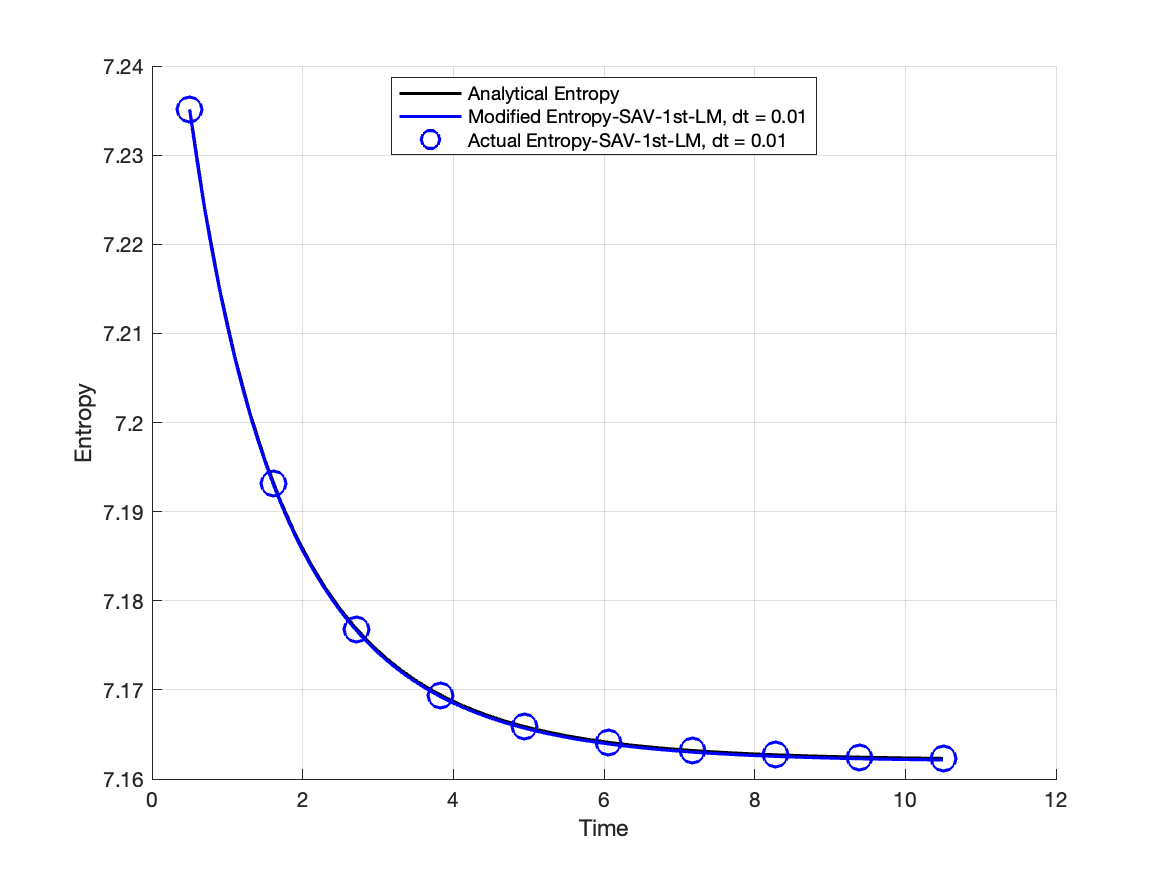}
        \caption{Entropy evolution with \( \Delta t = 0.01 \)}
    \end{subfigure}

     
    \begin{subfigure}[b]{0.45\textwidth}
        \includegraphics[width=\textwidth]{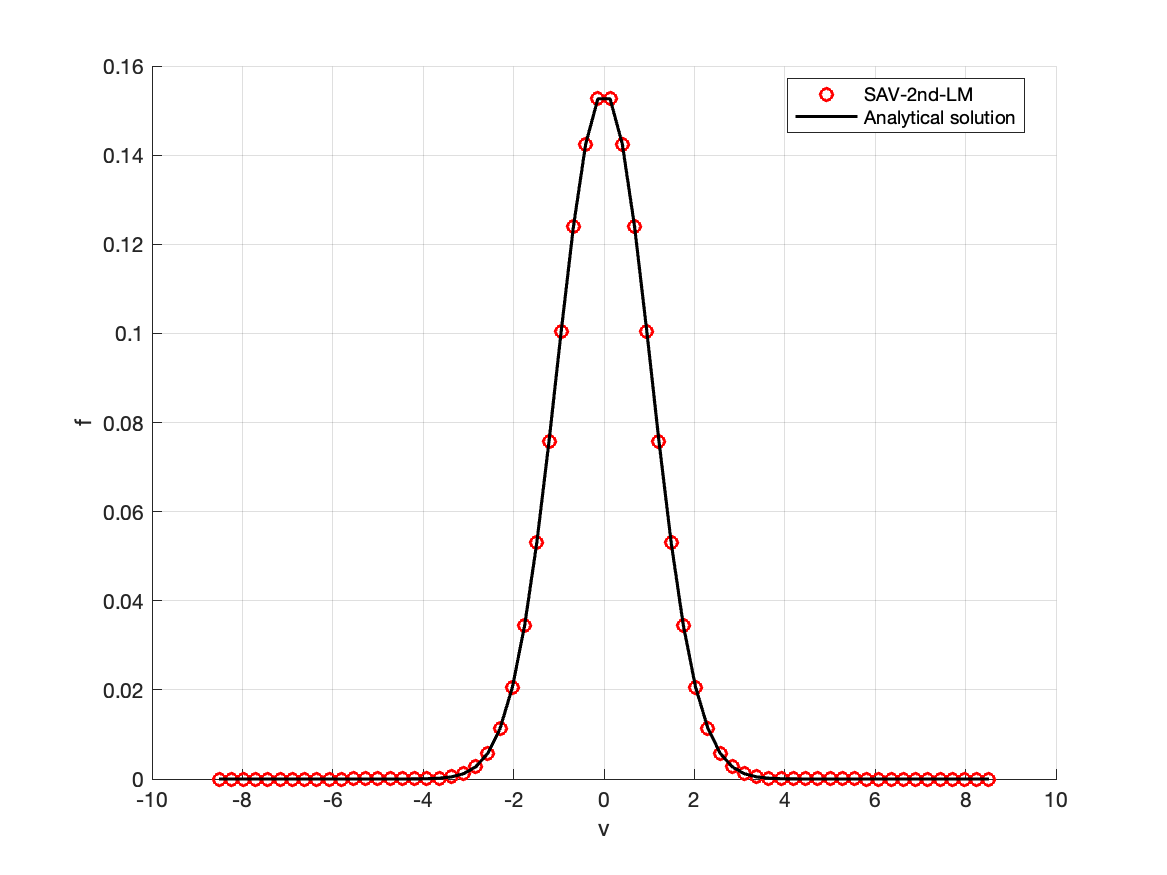}
        \caption{$f$ at $t_{\text{end}}=10.5$ with \( \Delta t = 0.2 \)}
    \end{subfigure}
    \hfill
    \begin{subfigure}[b]{0.45\textwidth}
        \includegraphics[width=\textwidth]{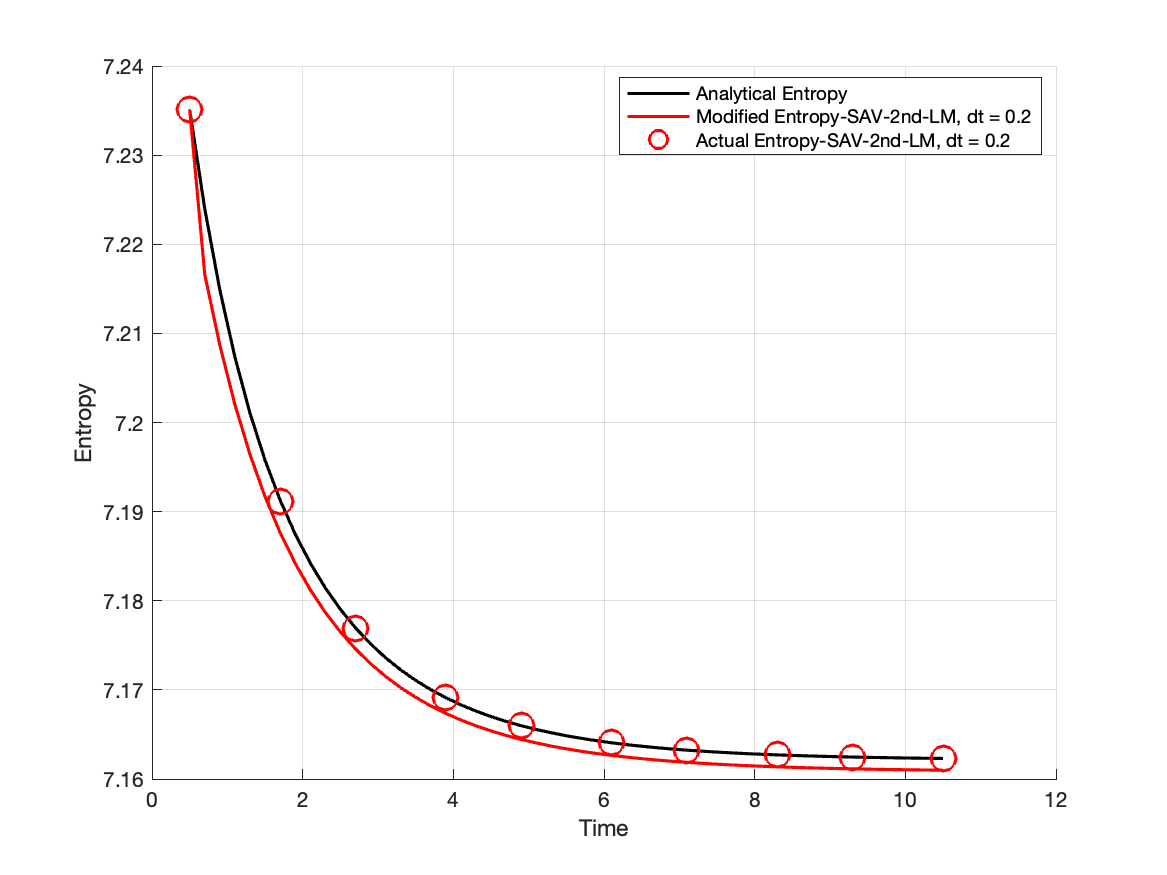}
        \caption{Entropy evolution with \( \Delta t = 0.2 \)}
    \end{subfigure}
    \hfill
    \begin{subfigure}[b]{0.45\textwidth}
        \includegraphics[width=\textwidth]{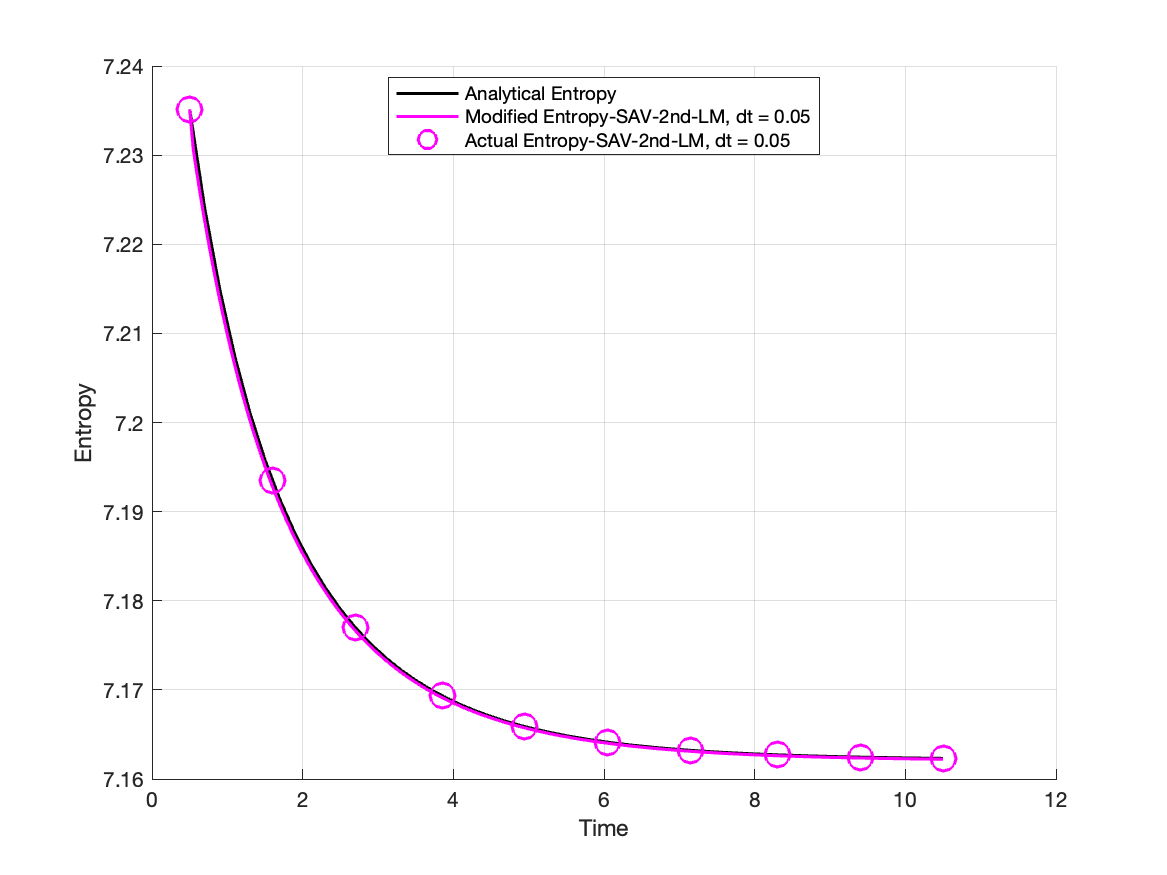}
        \caption{Entropy evolution with \( \Delta t = 0.05 \)}
    \end{subfigure}
    \hfill
    \begin{subfigure}[b]{0.45\textwidth}
        \includegraphics[width=\textwidth]{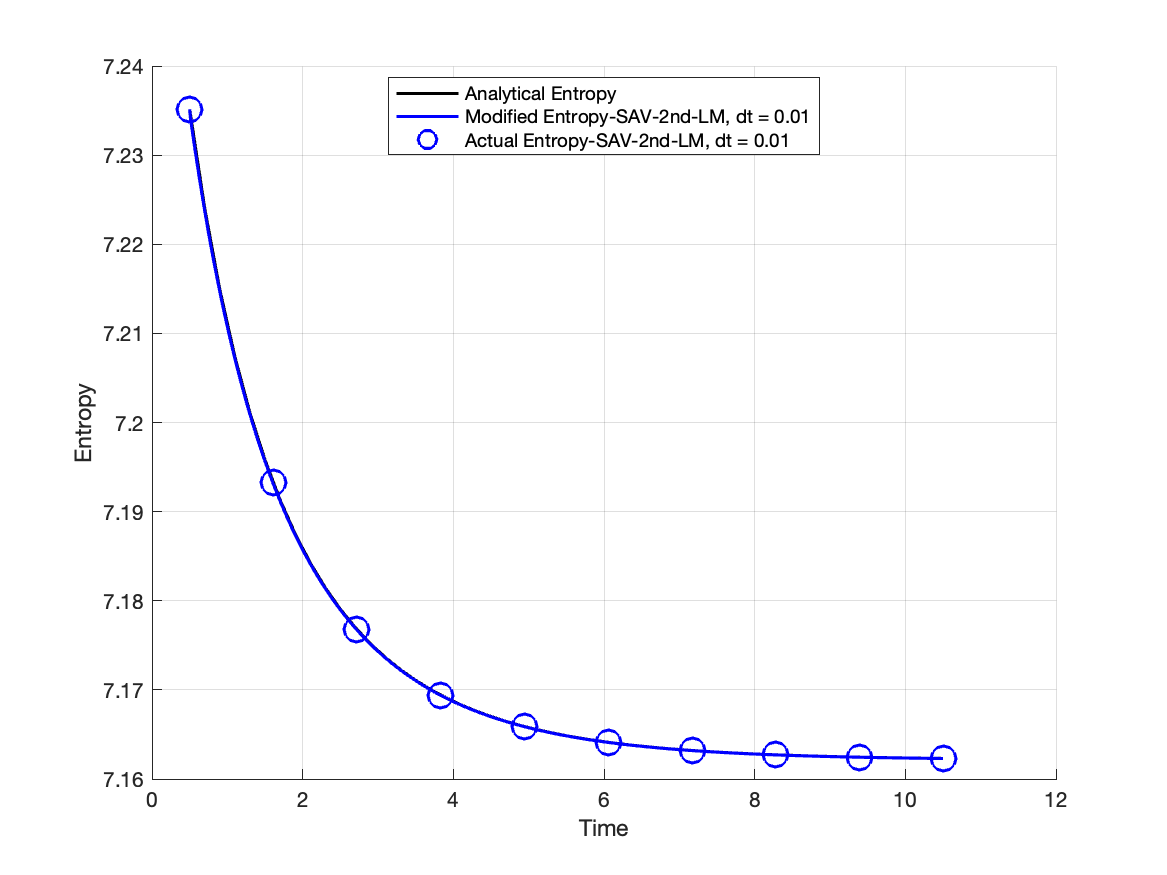}
        \caption{Entropy evolution with \( \Delta t = 0.01 \)}
    \end{subfigure}
\caption{Solution profile and entropy evolution for the Boltzmann equation with different time step sizes. Top 4 figures: \texttt{SAV-1st-LM}. Bottom 4 figures: \texttt{SAV-2nd-LM}.}
    \label{fig:boltz_comparison}
\end{figure}

Figure \ref{fig:boltz_comparison} showcases the results for the Boltzmann equation. The top 4 figures are obtained using \texttt{SAV-1st-LM}, while the bottom 4 figures are obtained using \texttt{SAV-2nd-LM}. These figures demonstrate that both methods accurately predict the solution, with the actual entropy closely matching the analytical entropy across all time steps. For larger time step sizes, the modified entropy exhibits a faster decay compared to the actual entropy, but the actual entropy consistently aligns closely with the entropy of the analytical solution. In fact, for all the results presented, the positivity and mass conservation correction is never triggered so the results are equivalent to those obtained by \texttt{SAV-1st} and \texttt{SAV-2nd}.

\begin{figure}[htbp!]
    \centering
    \begin{subfigure}[b]{0.45\textwidth}
        \includegraphics[width=\textwidth]{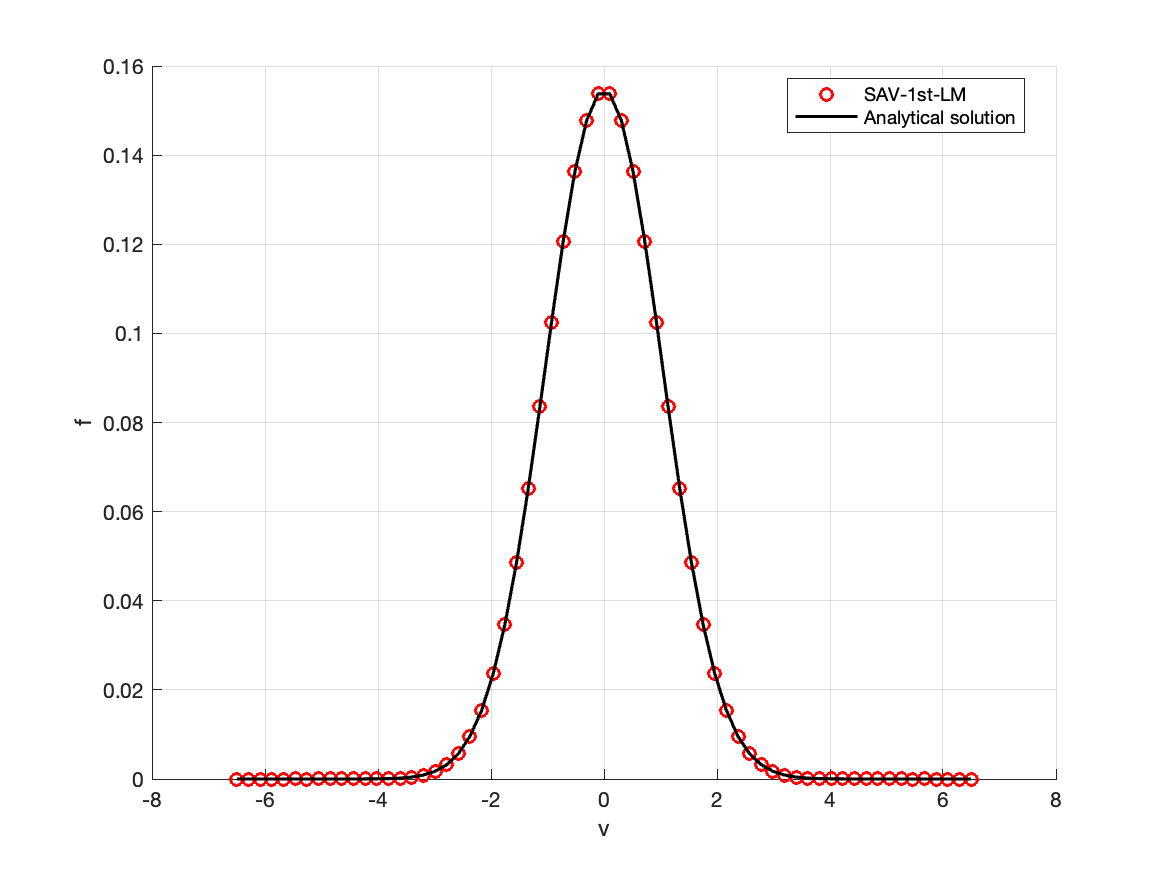}
        \caption{$f$ at $t_{\text{end}}=10.5$ with \( \Delta t  = 0.0041\)}
    \end{subfigure}
    \hfill
    \begin{subfigure}[b]{0.45\textwidth}
        \includegraphics[width=\textwidth]{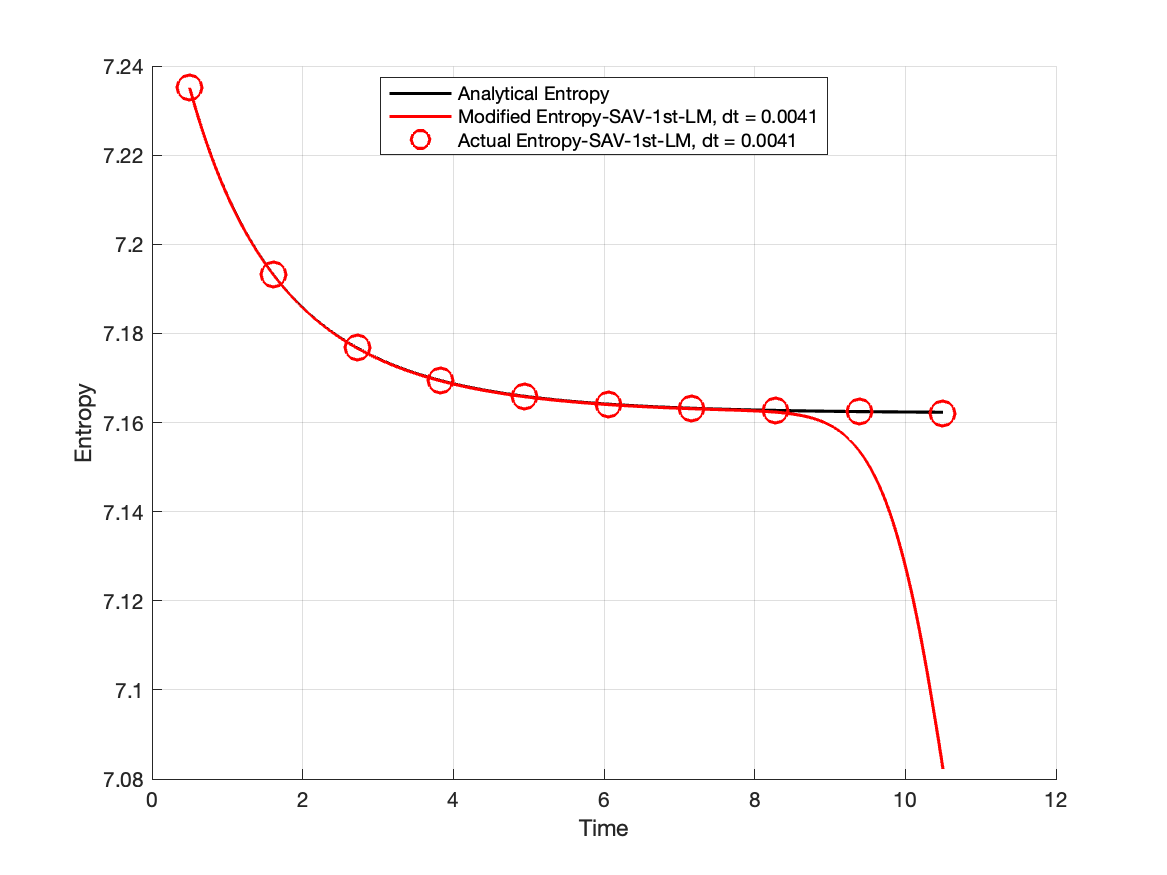}
        \caption{Entropy evolution with \( \Delta t = 0.0041 \)}
    \end{subfigure}
    \hfill
    \begin{subfigure}[b]{0.45\textwidth}
        \includegraphics[width=\textwidth]{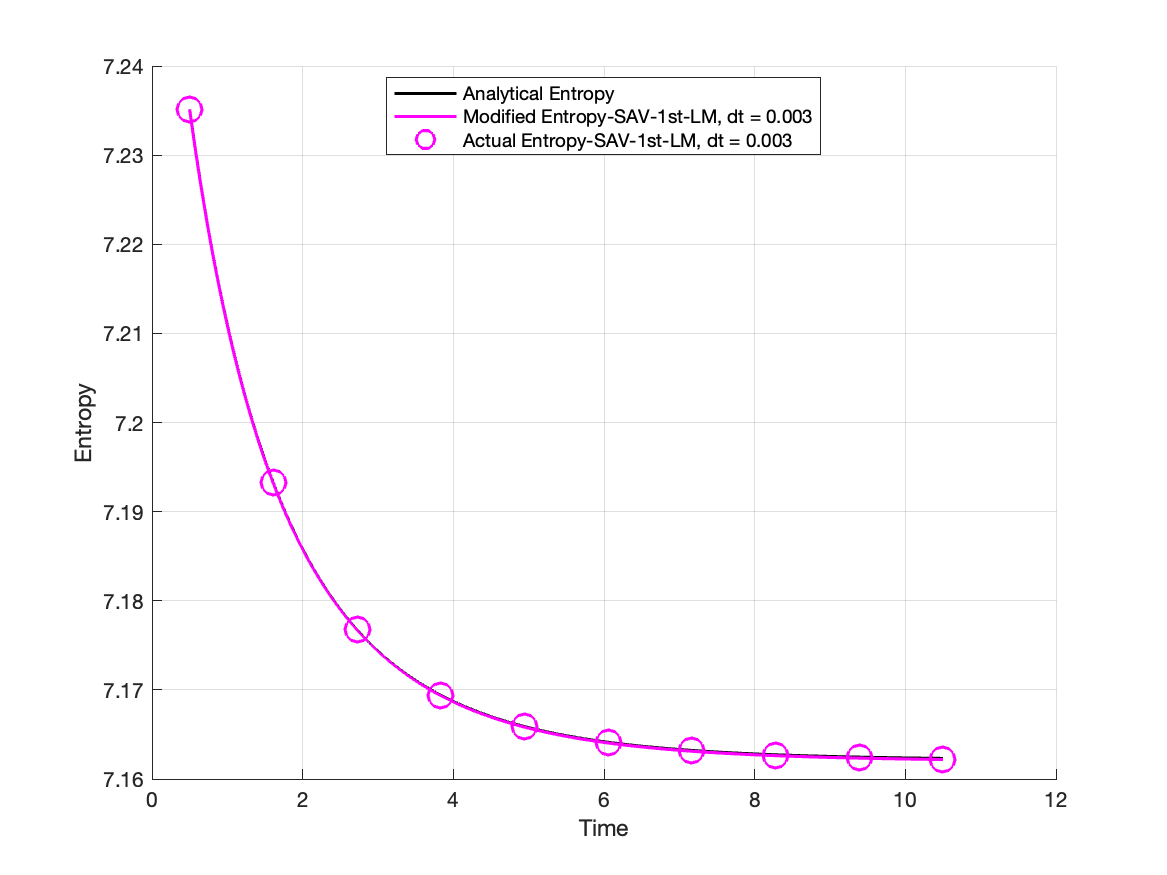}
        \caption{Entropy evolution with \( \Delta t = 0.003 \)}
    \end{subfigure}
    \hfill
    \begin{subfigure}[b]{0.45\textwidth}
        \includegraphics[width=\textwidth]{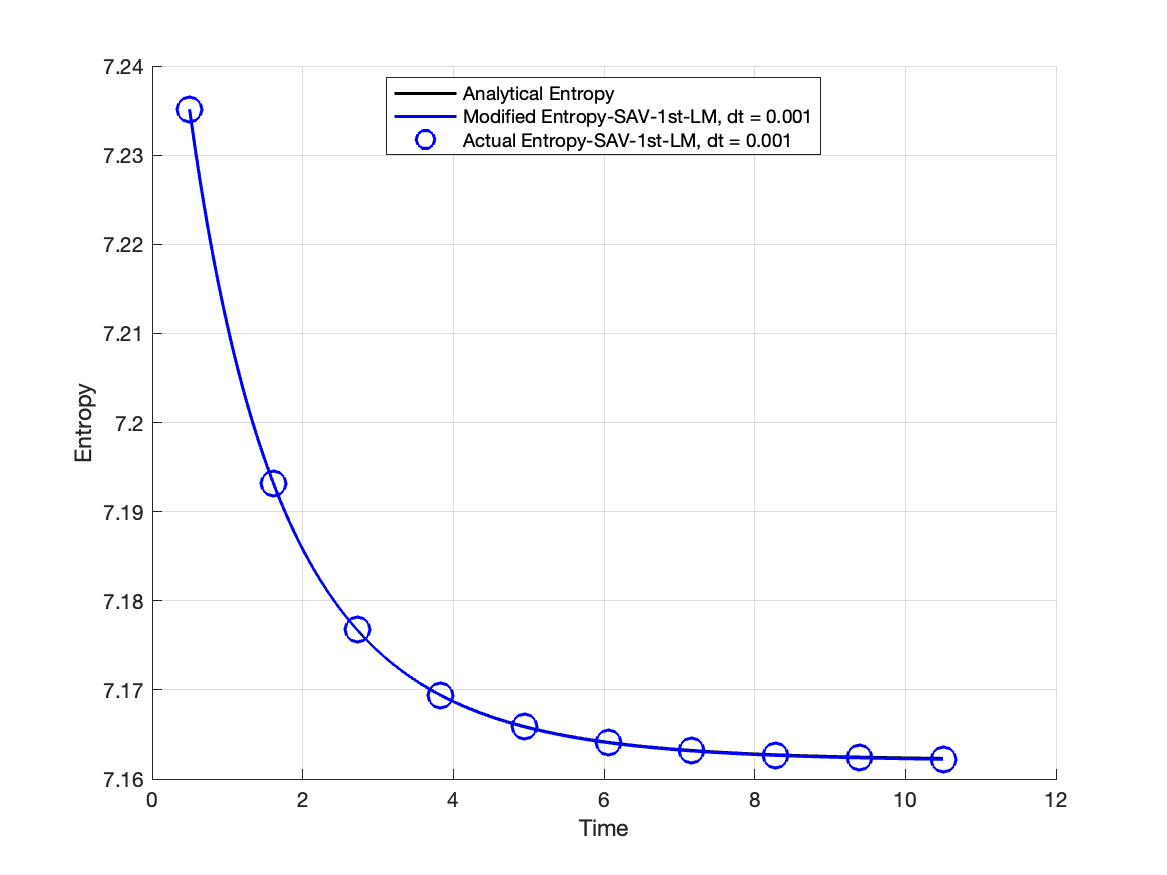}
        \caption{Entropy evolution with \( \Delta t = 0.001 \)}
    \end{subfigure}

    \begin{subfigure}[b]{0.45\textwidth}
        \includegraphics[width=\textwidth]{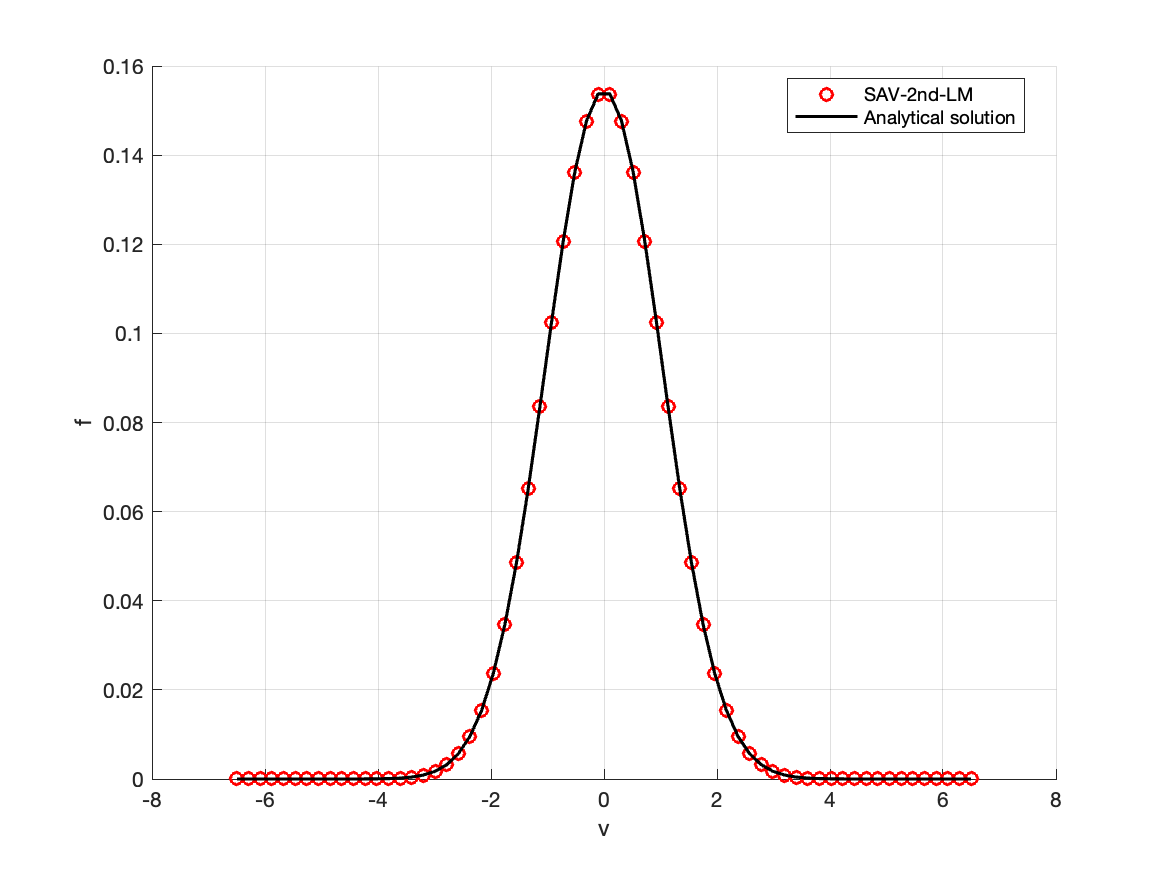}
        \caption{$f$ at $t_{\text{end}}=10.5$ with \( \Delta t  = 0.0035\)}
    \end{subfigure}
    \hfill
    \begin{subfigure}[b]{0.45\textwidth}
        \includegraphics[width=\textwidth]{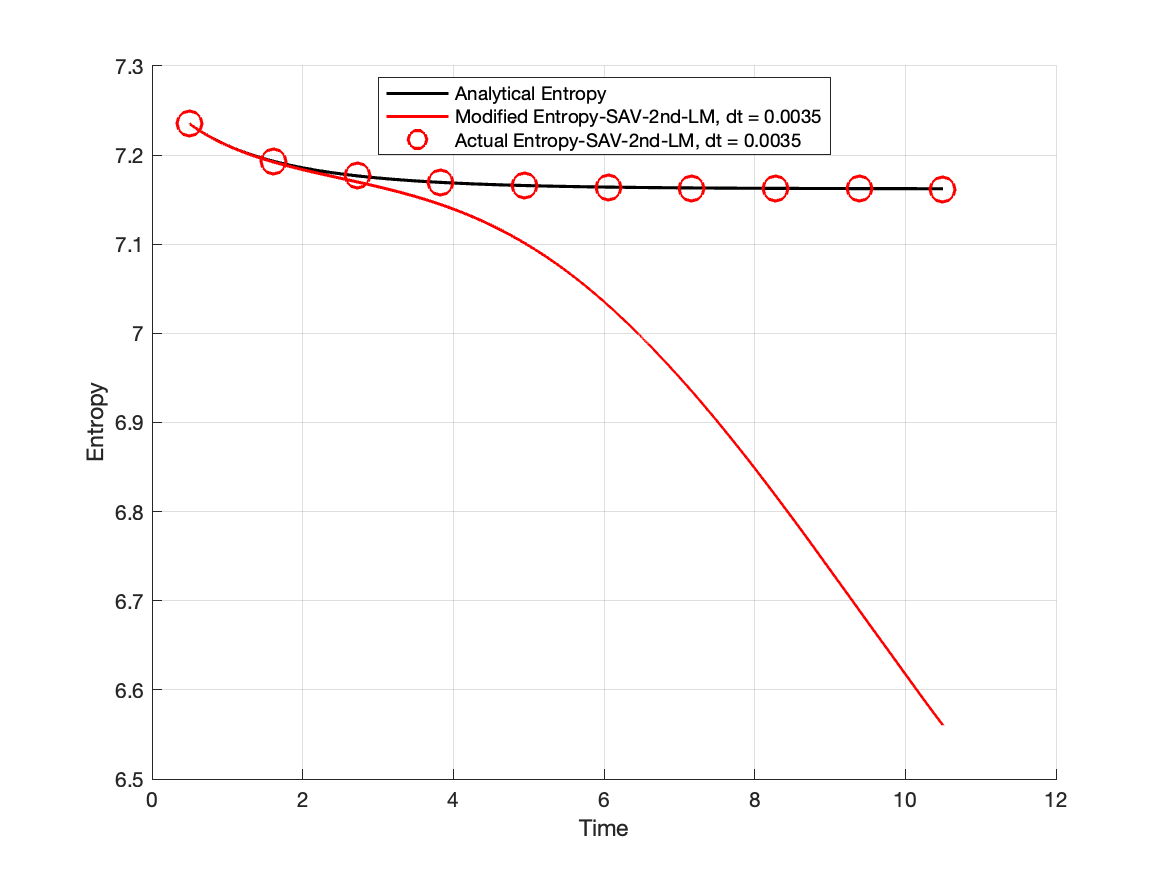}
        \caption{Entropy evolution with \( \Delta t = 0.0035 \)}
    \end{subfigure}
    \hfill
    \begin{subfigure}[b]{0.45\textwidth}
        \includegraphics[width=\textwidth]{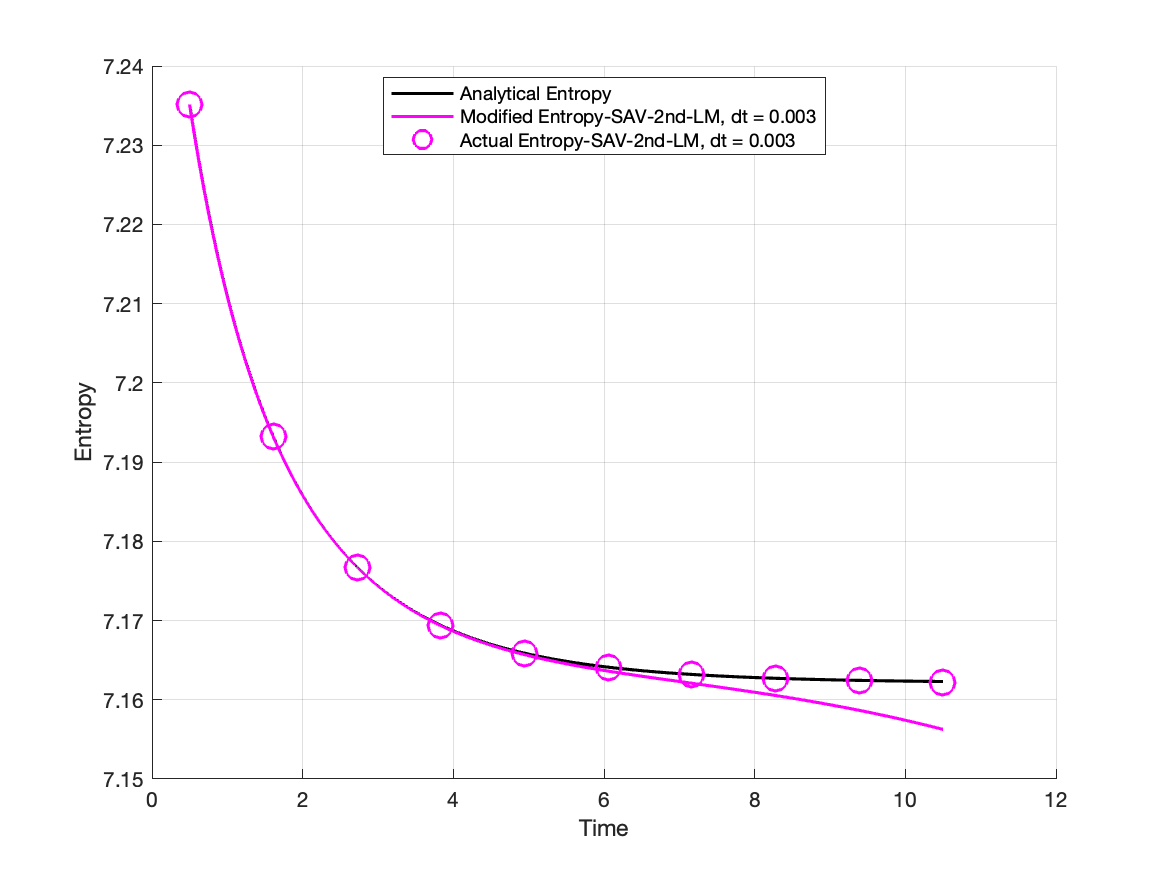}
        \caption{Entropy evolution with \( \Delta t = 0.003 \)}
    \end{subfigure}
    \hfill
    \begin{subfigure}[b]{0.45\textwidth}
        \includegraphics[width=\textwidth]{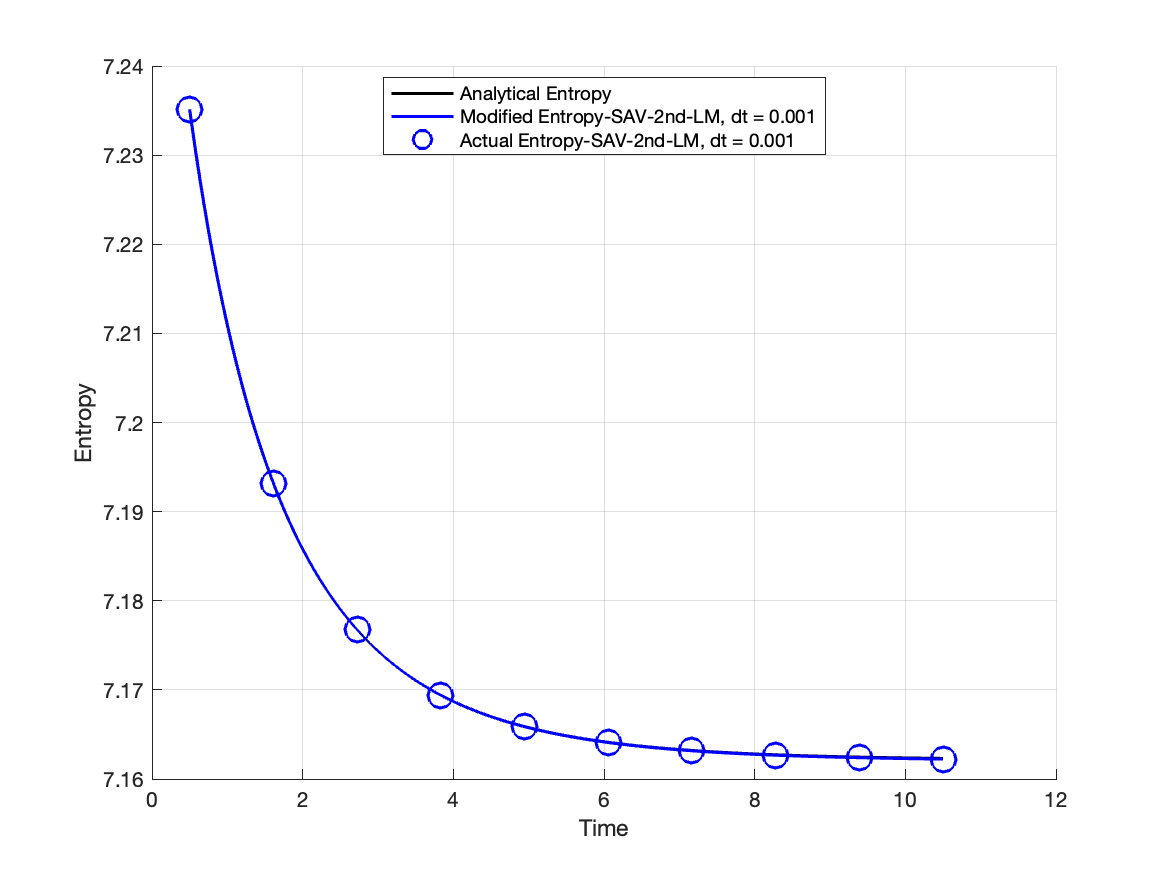}
        \caption{Entropy evolution with \( \Delta t = 0.001 \)}
    \end{subfigure}

    \caption{Solution profile and entropy evolution for the Landau equation with different time step sizes. Top 4 figures: \texttt{SAV-1st-LM}. Bottom 4 figures: \texttt{SAV-2nd-LM}.}
    \label{fig:landau_comparison}
\end{figure}

Figure \ref{fig:landau_comparison} showcases the results for the Landau equation. The top 4 figures are obtained using \texttt{SAV-1st-LM}, while the bottom 4 figures are obtained using \texttt{SAV-2nd-LM}. It is noteworthy that whenever the modified entropy deviates from the actual entropy, it often indicates that the positivity and mass conservation correction is triggered, in which case \texttt{SAV-1st} and \texttt{SAV-2nd} will fail. For larger time step sizes,  the modified entropy decays faster than the actual entropy, but the actual entropy always matches closely the analytical entropy.


Generally speaking, when using the SAV schemes with a sufficiently small time step size, the modified entropy will closely match the actual entropy. For larger time step sizes, the modified entropy is likely to decay at a faster rate than the actual entropy. This accelerated decay of the modified entropy contributes to stabilizing the scheme, enabling the actual entropy generated by the SAV schemes to closely approximate the true entropy.

\subsubsection{Positivity-preserving SAV scheme for the Boltzmann equation}

We now examine the performance of the first order positivity-preserving SAV scheme, \texttt{SAV-1st-P-B}, for solving the Boltzmann equation.

\begin{figure}[htbp!]
\centering
\begin{subfigure}{0.7\textwidth}
\centering
    \includegraphics[width=\linewidth]{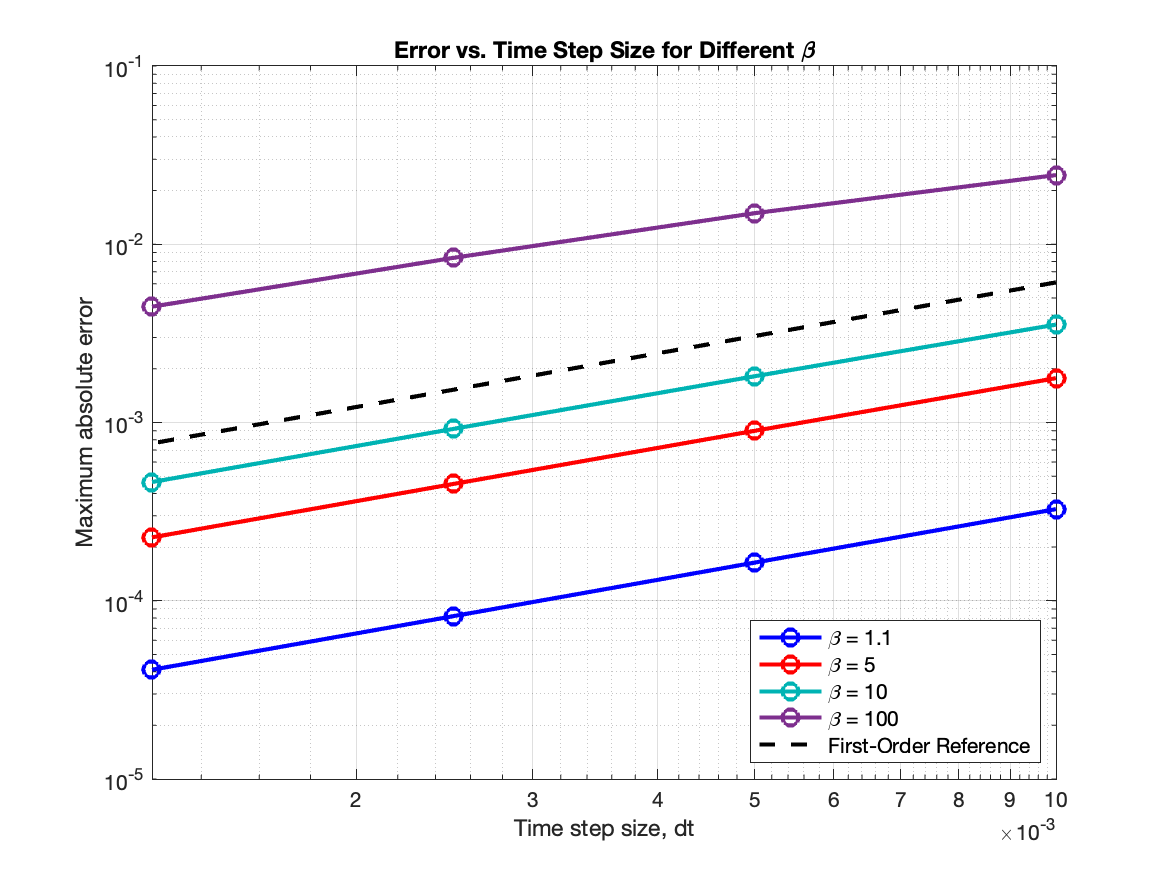}
    \caption{Convergence test for \texttt{SAV-1st-P-B} with different values of \(\beta\) (\(\beta = 1.1, 5, 10, 100\)) and \(\Delta t = 0.2, 0.1, 0.05, 0.025\).}
    \label{fig:error_vs_dt_psav1_boltz_lambda_comparison}
\end{subfigure}
\hfill
\begin{subfigure}{0.7\textwidth}
\includegraphics[width=\textwidth]{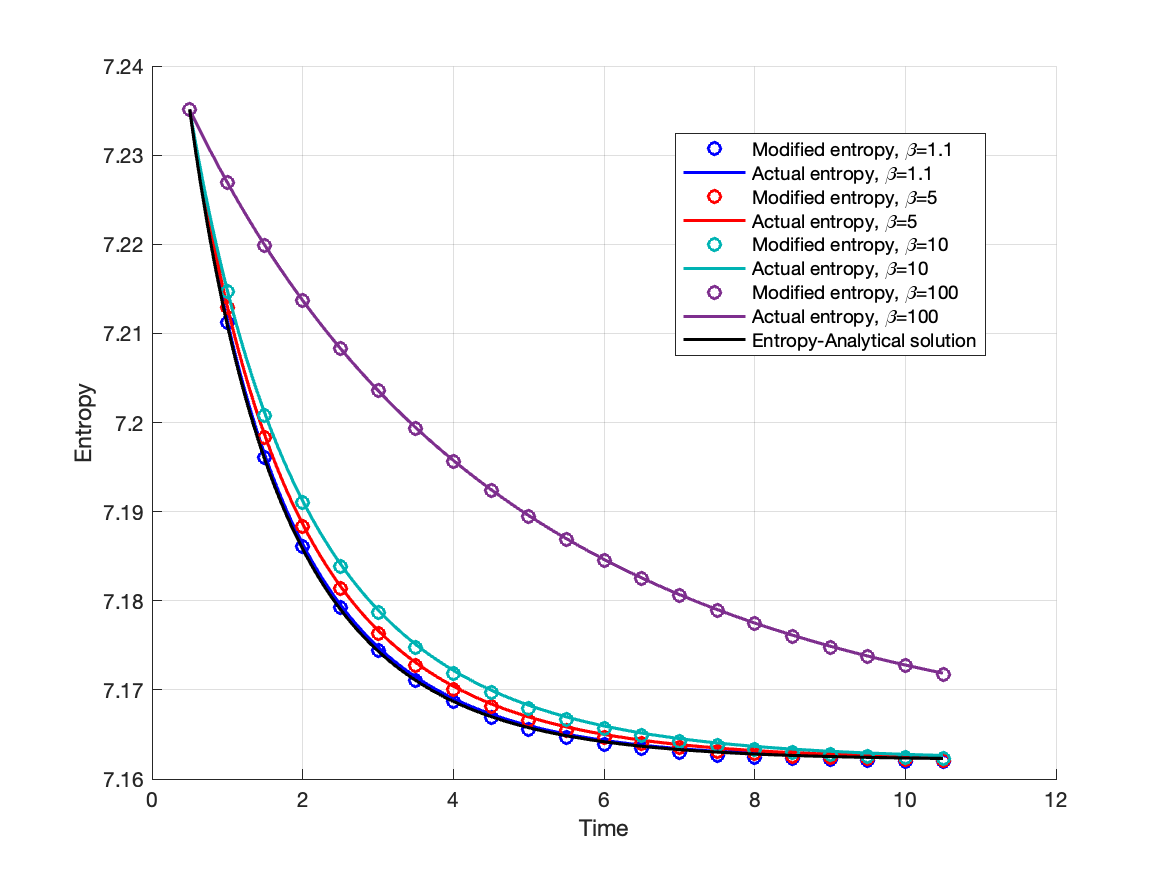}
\caption{Comparison of modified entropy, actual entropy, and entropy of the analytical solution for \texttt{SAV-1st-P-B} with different values of \(\beta\) (\(\beta = 1.1, 5, 10, 100\)) and a fixed time step size \(\Delta t = 0.025\).}
\label{fig:sav_entropy_lambda}
\end{subfigure}
\hfill
\caption{Convergence test and entropy evolution for \texttt{SAV-1st-P-B}.}
\label{fig:SAV-1st-P-B}
\end{figure}

For this scheme to be positive, the parameter $\beta$ is required to be larger than $\frac{r^0}{\sqrt{H_{\min}}}\max Q_B^-(f)$, where  \(r^0 = 2.6898\), \(\sqrt{H_{\min}} \approx 2.6780\), and \(\max Q_B^-(f) = 1\) in the BKW test. Therefore, \( \beta \) is chosen as 1.1, 5, 10, and 100. We first perform a convergence test using time step sizes of \( \Delta t = 0.2, 0.1, 0.05, \text{and } 0.025 \). For each time step size, we run the solution to $t_{\text{end}}=2.5$ and evaluate the maximum norm of the error between the numerical solution and the analytical one. The results are shown in Figure~\ref{fig:error_vs_dt_psav1_boltz_lambda_comparison}, from which we observe the expected first order convergence. It is also clear that larger $\beta$ results in larger error in magnitude. So in practice, $\beta$ should be chosen as close as possible to the required lower bound. We then fix the time step size \(\Delta t = 0.025\) and plot the entropy evolution for different values of $\beta$ in Figure~\ref{fig:sav_entropy_lambda}. The modified entropy always matches well the actual entropy. However, larger $\beta$ results in larger deviation from the analytical entropy so the proper choice of $\beta$ is also critical to obtain a correct entropy dissipation rate.



\subsection{Test case 2:}

In this test, we consider the following initial condition 
\begin{equation}\label{initial ex2}
f^0(v)=\frac{\rho_1}{2\pi T_1}\exp\left(-\frac{|v-V_1|^2}{2T_1}\right)+\frac{\rho_2}{2\pi T_2}\exp\left(-\frac{|v-V_2|^2}{2T_2}\right),
\end{equation}
with $\rho_1=\rho_2=1/2$, $T_1=T_2=1$, and $V_1=(-1,2)$, $V_2=(3,-3)$. 
For the Boltzmann equation, $S$ is set to 5; and for the Landau equation, $S$ is set to 7.5. For this initial condition, we don't know the analytical solution but we do know that the solution will relax to the Maxwellian after a long time.

\begin{figure}[htbp]
\centering
\begin{subfigure}{0.48\textwidth}
    \includegraphics[width=\linewidth]{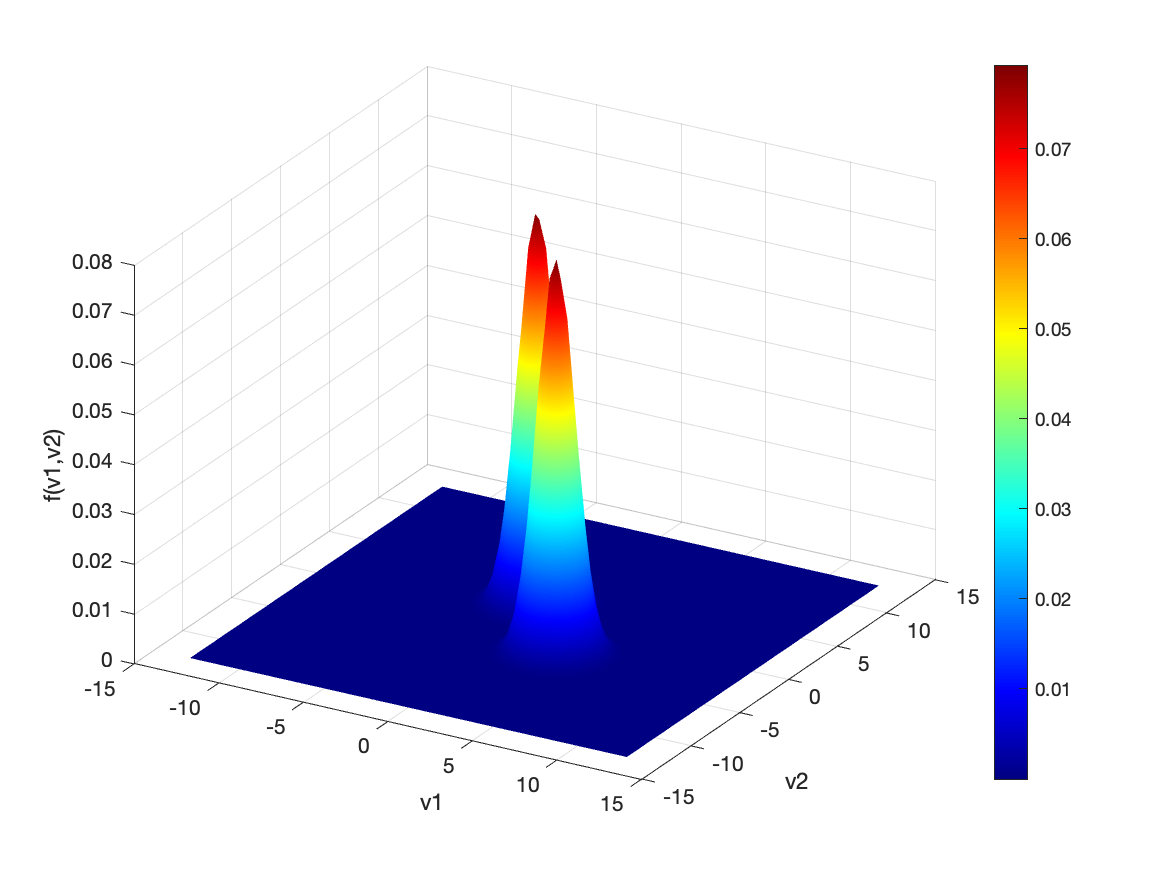}
    \caption{Boltzmann equation: solution at $t_0=0$}
    \label{fig:exm2_f0}
\end{subfigure}
\hfill
\begin{subfigure}{0.48\textwidth}
    \includegraphics[width=\linewidth]{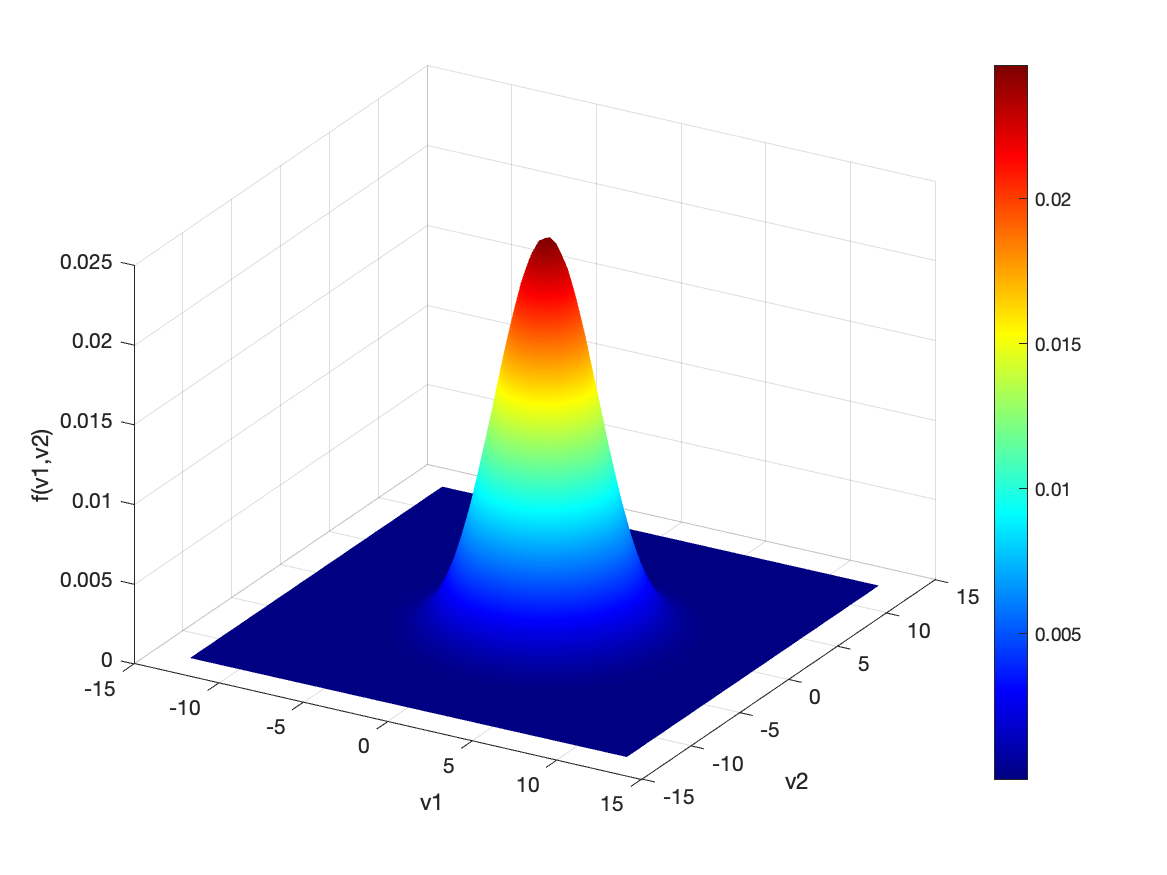}
    \caption{Boltzmann equation: solution at $t=10$}
    \label{fig:exm2_ft}
\end{subfigure}
\hfill
\begin{subfigure}{0.48\textwidth}
    \includegraphics[width=\linewidth]{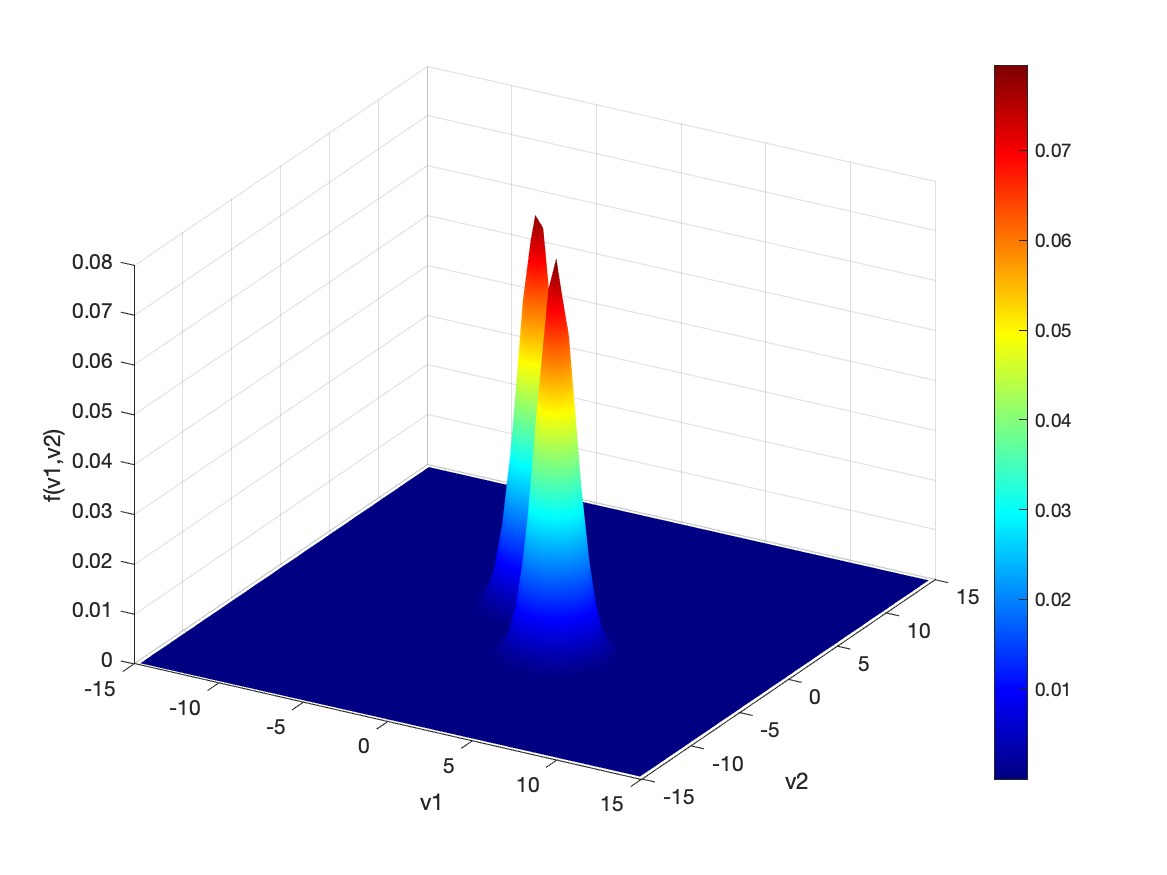}
    \caption{Landau equation: solution at $t_0=0$}
    \label{fig:f0-landau}
\end{subfigure}
\hfill
\begin{subfigure}{0.48\textwidth}
    \includegraphics[width=\linewidth]{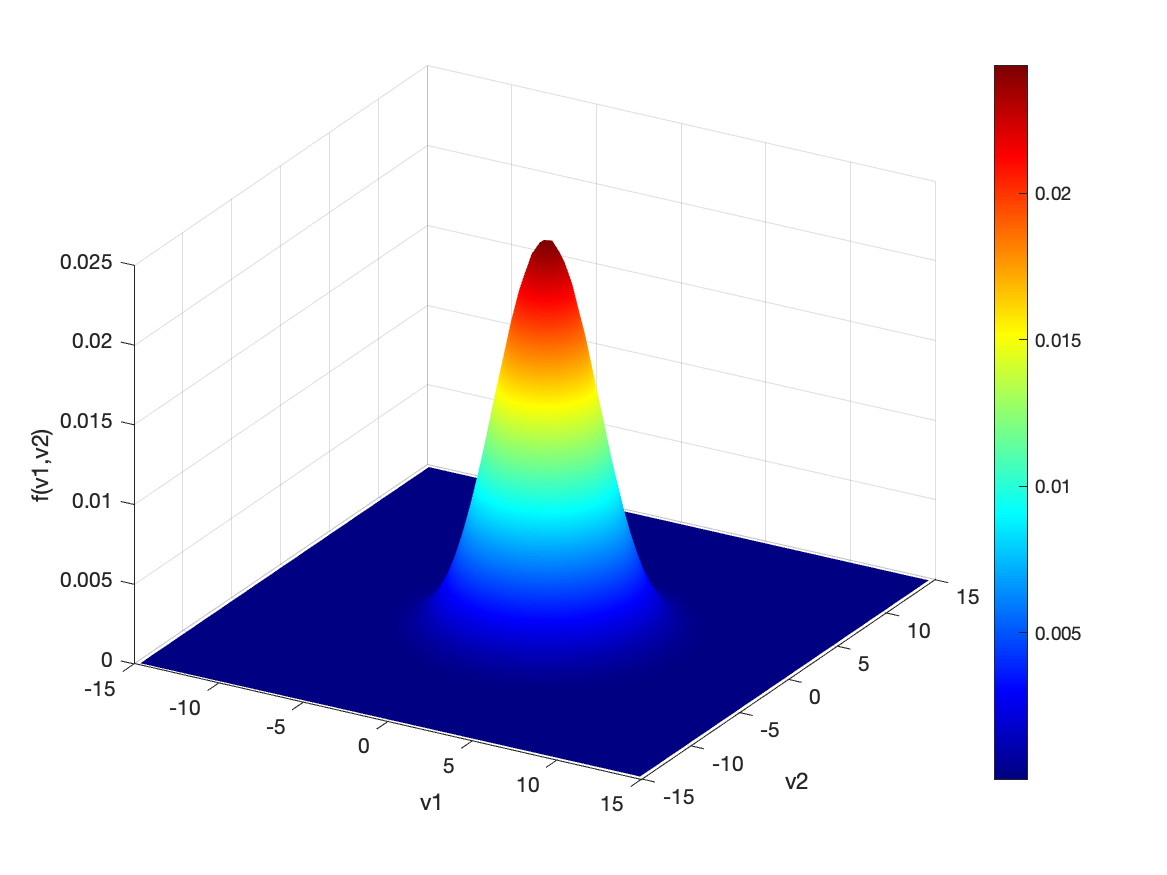}
    \caption{Landau equation: solution at $t=10$}
    \label{fig:ft-landau}
\end{subfigure}
\caption{Solutions computed using \texttt{SAV-2nd-LM} for the Boltzmann equation (with $\Delta t=0.01$) and Landau equation (with $\Delta t=0.002$).}
\label{fig:exm2-1}
\end{figure}

\begin{figure}[htbp]
\centering
\begin{subfigure}{0.48\textwidth}
    \includegraphics[width=\linewidth]{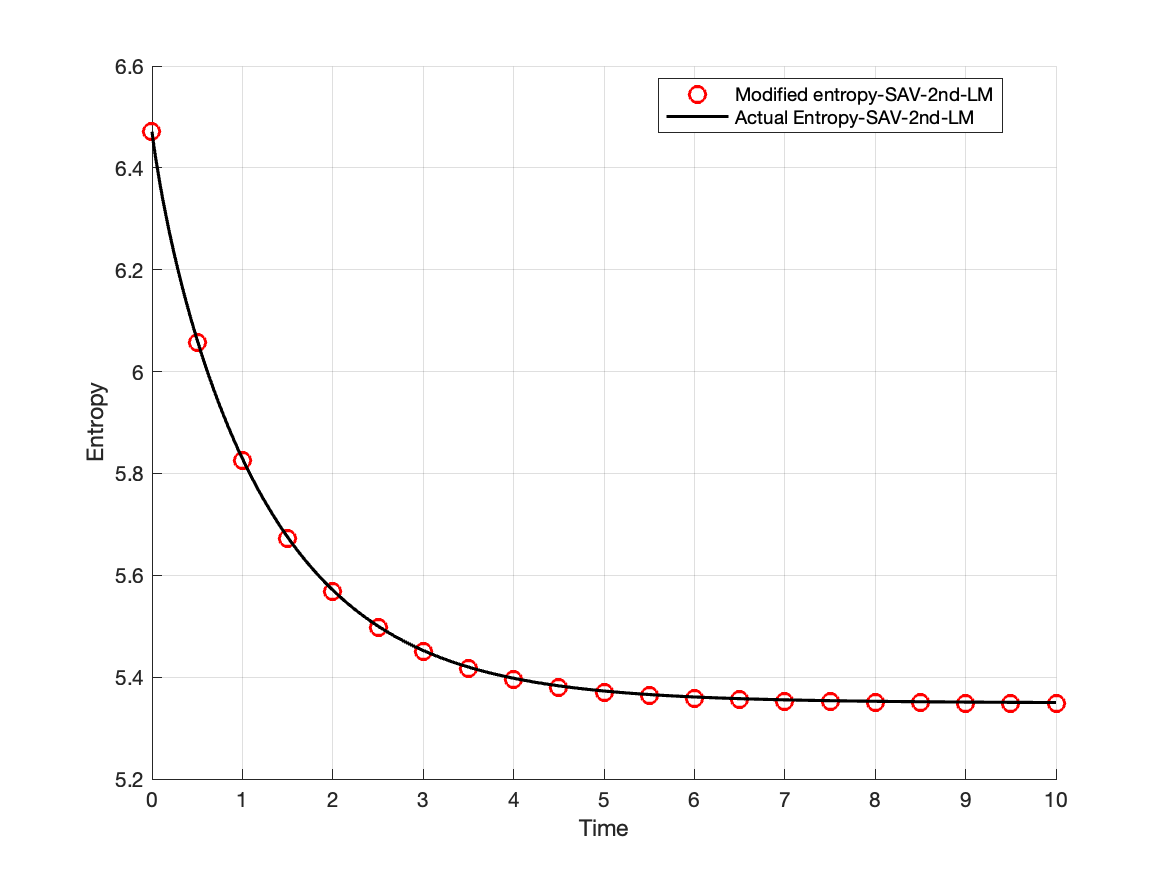}
    \caption{Boltzmann equation: entropy decay.}
    \label{fig:exm2_fsav2euler}
\end{subfigure}
\hfill
\begin{subfigure}{0.48\textwidth}
    \includegraphics[width=\linewidth]{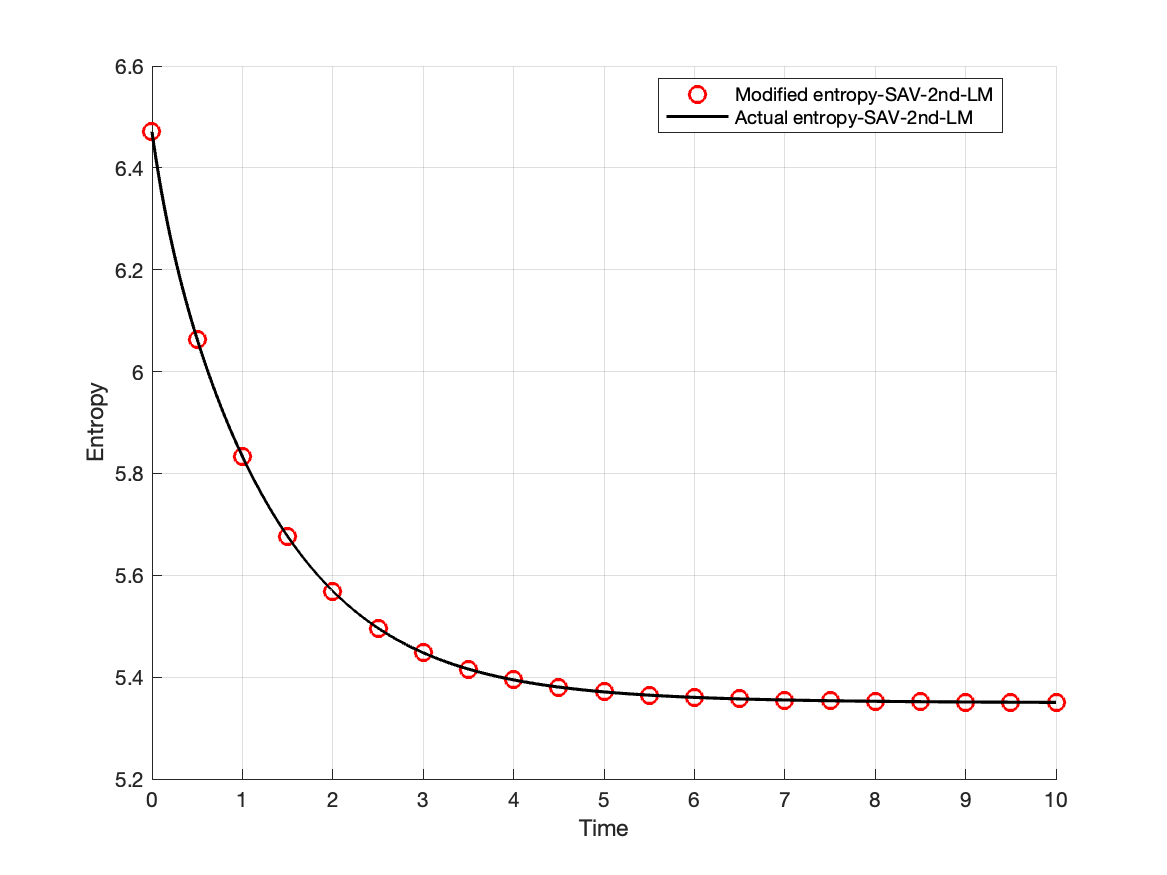}
    \caption{Landau equation: entropy decay.}
    \label{fig:exm2-entropy-landau}
\end{subfigure}
\caption{Entropy evolution for the Boltzmann equation (with $\Delta t=0.01$) and Landau equation (with $\Delta t=0.002$) computed using \texttt{SAV-2nd-LM}.}
\label{fig:exm2-3}
\end{figure}

In  Figure~\ref{fig:exm2-1}, we illustrate the evolution of the distribution function for both the Boltzmann equation (with $\Delta t=0.01$) and Landau equation (with $\Delta t=0.002$) until final time $t=10$, computed using \texttt{SAV-2nd-LM}. The expected trend of the solution is observed. In Figure~\ref{fig:exm2-3}, we compared the modified entropy against the actual entropy obtained. The results demonstrate a consistent alignment between the modified and actual entropy, highlighting the ability of the proposed scheme to capture the entropy dissipation structure.

\section{Conclusions}
\label{sec:con}

We developed novel numerical schemes to tackle the dual challenge of enabling entropy dissipation and preserving positivity for general kinetic equations, leveraging the recently introduced SAV approach. Both the first order and second order schemes were constructed. We applied the proposed schemes to the nonlinear Boltzmann equation and Landau equation, which are among the most challenging kinetic equations, and presented convincing numerical results which showed that the proposed schemes are both robust and efficient. Future work includes the extension of these schemes to treat the spatially inhomogeneous kinetic equations.




\bibliographystyle{siamplain}
\bibliography{hu_bibtex}
\end{document}